\documentclass[a4paper,12pt]{article}
\usepackage[utf8]{inputenc}
\usepackage[T1]{fontenc}
\usepackage{amsmath}
\usepackage{emptypage} 
\usepackage{a4wide}
\usepackage[pagestyles]{titlesec}
\usepackage{accents}
\newcommand{\ubar}[1]{\underaccent{\bar}{#1}}
\newpagestyle{main}{
\headrule
\sethead[\thepage][\chaptertitle][] 
{}{\sectiontitle}{\thepage} 
}

\usepackage{latexsym}
\usepackage{amssymb}
\usepackage{amsthm}
\usepackage{graphicx}
\usepackage{amstext}
\usepackage{bbm}
\usepackage{xcolor}
\usepackage{hyperref}
\usepackage[capitalize
]{cleveref}
\usepackage[math]{cellspace}

\makeatletter
\def\underbracex#1#2{\mathop{\vtop{\m@th\ialign{##\crcr
   $\hfil\displaystyle{#2}\hfil$\crcr
   \noalign{\kern3\p@\nointerlineskip}%
   #1\crcr\noalign{\kern3\p@}}}}\limits}

\def\upbracefilla{$\m@th \setbox\z@\hbox{$\braceld$}%
  \bracelu\leaders\vrule \@height\ht\z@ \@depth\z@\hfill 
\kern\p@\vrule \@width\p@\kern\p@\vrule \@width\p@\kern\p@\vrule \@width\p@
$}

\def\upbracefillb{$\m@th \setbox\z@\hbox{$\braceld$}%
\vrule \@width\p@\kern\p@\vrule \@width\p@\kern\p@\vrule \@width\p@\kern\p@
 \leaders\vrule \@height\ht\z@ \@depth\z@\hfill\bracerd
  \braceld\leaders\vrule \@height\ht\z@ \@depth\z@\hfill
\kern\p@\vrule \@width\p@\kern\p@\vrule \@width\p@\kern\p@\vrule \@width\p@
$}

\def\upbracefillc{$\m@th \setbox\z@\hbox{$\braceld$}%
\vrule \@width\p@\kern\p@\vrule \@width\p@\kern\p@\vrule \@width\p@\kern\p@
\leaders\vrule \@height\ht\z@ \@depth\z@\hfill
\kern\p@\vrule \@width\p@\kern\p@\vrule \@width\p@\kern\p@\vrule \@width\p@
$}

\def\upbracefill{$\m@th \setbox\z@\hbox{$\braceld$}%
\vrule \@width\p@\kern\p@\vrule \@width\p@\kern\p@\vrule \@width\p@\kern\p@
 \leaders\vrule \@height\ht\z@ \@depth\z@\hfill\braceru$}

\def\upbracefillL{$\m@th \setbox\z@\hbox{$\braceld$}%
  \bracelu\leaders\vrule \@height\ht\z@ \@depth\z@\hfill 
\kern\p@\vrule \@width\p@\kern\p@\vrule \@width\p@\kern\p@\vrule \@width\p@
$}

\def\upbracefillR{$\m@th \setbox\z@\hbox{$\braceld$}%
\vrule \@width\p@\kern\p@\vrule \@width\p@\kern\p@\vrule \@width\p@\kern\p@
 \leaders\vrule \@height\ht\z@ \@depth\z@\hfill\bracerd
  \braceld\leaders\vrule \@height\ht\z@ \@depth\z@\hfill
 \leaders\vrule \@height\ht\z@ \@depth\z@\hfill\braceru$}

\def\BRR{\underbracex\upbracefillLR}

\def\upbracefillLR{$\m@th \setbox\z@\hbox{$\braceld$}%
  \bracelu\leaders\vrule \@height\ht\z@ \@depth\z@\hfill 
\leaders\vrule \@height\ht\z@ \@depth\z@\hfill\bracerd
  \braceld\leaders\vrule \@height\ht\z@ \@depth\z@\hfill
 \leaders\vrule \@height\ht\z@ \@depth\z@\hfill\braceru
$}
\makeatother

\def\R{\mathbb R}

\def\N{\mathbb N}
\def\Z{\mathbb Z}

\newcommand{\p}{\psi}
\renewcommand{\t}{\theta}

\newcommand{\ti}{\tilde }

\newcommand{\tip}{{\tilde \p}}
\newcommand{\tit}{{\tilde \t}}

\newcommand{\hatp}{{\hat\p}}
\newcommand{\hatt}{{\hat \t}}

\newcommand{\tixi}{\ti \xi}
\newcommand{\tiu}{\ti \u}

\newcommand{\hxi}{\hat \xi}
\newcommand{\hu}{\hat \u}

\newcommand{\luen}{{\lambda_{n}^{\eps}}}

\newcommand{\tttn}{\theta_n^{\eps}(\xi(t),\u(t),x)}
\newcommand{\pttn}{\psi_n^{\eps}(\xi(t),\u(t),x)}
\newcommand{\ttn}{\theta_n^{\eps}(\xi,\u,x)}

\newcommand{\ttbn}{\theta_n^{\eps}(\bar\xi,\bar\u,x)}

\newcommand{\ptbn}{\p_n^{\eps}(\bar\xi,\bar\u,x)}

\newcommand{\ptn}{\p_n^{\eps}(\xi,\u,x)}

\def\u{u}
\def\g{\gamma}
\def\gz{\gamma(\u_0)}

\def\Z{\g(x-\xi)}

\def\a{\alpha}

\renewcommand{\v}{v}
\newcommand{\w}{w}

\newcommand{\pa}{\partial_}

\newcommand{\tp}{\l(\begin{matrix}
\t\\
\p\\
\end{matrix}\r)}

\newcommand{\tpP}{\l(\begin{matrix}
\t'\\
\p'\\
\end{matrix}\r)}

\newcommand{\dt}{\pa t}
\newcommand{\ddt}{\fr d{dt}}

\newcommand{\ds}{\fr d {ds}}

\newcommand{\dxi}{\pa \xi}

\newcommand{\du}{\pa \u}

\newcommand{\dx}{\pa x}
\newcommand{\dZ}{\pa Z}
\newcommand{\deps}{\pa \eps}
\def\xx{_{xx}}

\renewcommand{\d}{\dot}

\def \F{F(\eps,x)}

\def \Go{{\cal G}_1^{\eps}}
\def \Gt{{\cal G}_2^{\eps}}

\def \Gn{{\cal G}_n^{\eps}}

\def \Rn{{\cal R}_n^{\eps}}

\def \Rtno{[{\cal R}_n^{\eps}(\xi,\u,x)]_1}

\def \Rtnt{[{\cal R}_n^{\eps}(\xi,\u,x)]_2}

\def \L{{\cal L}_{\xi,\u}}

\def \hatL{ {\hat{\cal L}}_{\xi,\u} }

\def \tiGo{\tilde{\cal G}_1^{\eps}}
\def \tiGt{\tilde{\cal G}_2^{\eps}}

\def \tiGn{\tilde{\cal G}_n^{\eps}}
\DeclareMathOperator{\sech}{sech}

\newcommand{\BR}[2]{\BRR{{#1}}_\text{{#2}}}

\newcommand{\OM}[2]{\Omega \l({{#1}},{{#2}}\r)}
\newcommand{\OT}[2]{\stackrel{{\smash{\scriptscriptstyle{#2}}}}{{#1}}}
\newcommand{\Ltwo}[2]{\l\langle {#1},{#2} \r\rangle_{L^2(\R)\oplus L^2(\R)}}

\newcommand{\ltwortwoaxiZ}[2]{\l\langle {#1},{#2} \r\rangle_{L^{2,a}_{\xi,Z}(\R^2)}}

\newcommand{\ltwortwoaxix}[2]{\l\langle {#1},{#2} \r\rangle_{L^{2,a}_{\xi,x}(\R^2)}}

\newcommand{\Ltwortwoa}[2]{\l\langle {#1},{#2} \r\rangle_{L^{2,a}(\R^2)\oplus L^{2,a}(\R^2)}}

\newcommand{\Ltwortwoaxix}[2]{\l\langle {#1},{#2} \r\rangle_{L^{2,a}_{\xi,x}(\R^2)\oplus L^{2,a}_{\xi,x}(\R^2)}}

\newcommand{\ltwoa}[2]{\l\langle {#1},{#2} \r\rangle_{L^{2,a}(\R)}}

\newcommand{\tnv}[3]{\l|{#1}\r|_{L^\infty({[{#2},{#3}]},H^1(\R))}}
\newcommand{\tnvinf}[3]{\l|{#1}\r|_{L^\infty(\R)L^\infty({[{#2},{#3}]})}}

\newcommand{\tnw}[3]{\l|{#1}\r|_{L^\infty({[{#2},{#3}]},L^2(\R))}}
\newcommand{\nw}[1]{\l|{#1}\r|_{L^2(\R)}}

\newcommand{\nltwox}[1]{\l|{#1}\r|_{L^2_x(\R)}}

\newcommand{\nhonex}[1]{\l|{#1}\r|_{H^1_x(\R)}}

\newcommand{\nltwo}[1]{\l|{#1}\r|_{L^2(\R)}}
\newcommand{\nlinf}[1]{\l|{#1}\r|_{L^\infty(\R)}}

\newcommand{\nhtwortwoxix}[1]{\l|{#1}\r|_{H^{2}_{\xi,x}(\R^2)}}

\newcommand{\nhonertwoxix}[1]{\l|{#1}\r|_{H^1_{\xi,x}(\R^2)}}

\newcommand{\nltwortwoxiZ}[1]{\l|{#1}\r|_{L^{2}_{\xi,Z}(\R^2)}}
\newcommand{\nltwortwoaxiZ}[1]{\l|{#1}\r|_{L^{2,\a}_{\xi,Z}(\R^2)}}

\newcommand{\nhtwortwoxiZ}[1]{\l|{#1}\r|_{H^{2}_{\xi,Z}(\R^2)}}

\newcommand{\nltwortwoa}[1]{\l|{#1}\r|_{L^{2,\a}(\R^2)}}

\newcommand{\nhtwoa}[1]{\l|{#1}\r|_{H^{2,\a}(\R)}}

\newcommand{\nhone}[1]{\l|{#1}\r|_{H^1(\R)}}

\newcommand{\CR}{\color{black}}
\newcommand{\fr}{\frac}

\def\be#1\ee{\begin{align}#1\end{align}}
\def\ben#1\een{\begin{align+}#1\end{align+}}
\def\ba#1\ea{\begin{aligned}[t]#1\end{aligned}}
\def\bs#1\es{\begin{split}#1\end{split}}

\renewcommand{\l}{\left}
\renewcommand{\r}{\right}
\def\?{????}

\newcommand{\bma}{\begin{pmatrix}}
\newcommand{\ema}{\end{pmatrix}}

\newcommand{\bmat}{\begin{bmatrix}}
\newcommand{\emat}{\end{bmatrix}}

\newcommand{\bca}{\begin{cases}}
\newcommand{\eca}{\end{cases}}

\def\J{\mathbb J}
\newcommand{\re}{\eqref}
\newcommand{\nn}{\nonumber}

\newcommand{\la}{\label}
\newcommand{\beq}{\begin{equation}}
\newcommand{\eeq}{\end{equation}}
\newcommand{\eps}{\varepsilon}

\renewcommand{\qed}{\protect~\protect\hfill $\Box$}

\newcommand{\ran}{{\rm ran}\,}
\renewcommand{\ker}{{\rm ker}\,}

\begin{document}


\newtheoremstyle{plainNEW}
  {}
  {}
  {\itshape}
  {}
  {\boldmath\bfseries}
  {.}
  { }
  {\thmname{#1}\thmnumber{ #2}\thmnote{ (#3)}}%


\theoremstyle{plainNEW}

\newtheorem{theorem}{Theorem}[section]

\newtheorem{definition}[theorem]{Definition}
\newtheorem{deflem}[theorem]{Definition and Lemma}
\newtheorem{lemma}[theorem]{Lemma}
\newtheorem{corollary}[theorem]{Corollary}
\newcommand{\bde}{\begin{definition}}
\newcommand{\ede}{\end{definition}}
\newcommand{\ble}{\begin{lemma}}
\newcommand{\ele}{\end{lemma}}
\newcommand{\bre}{\begin{remark}}
\newcommand{\ere}{\end{remark}}
\newcommand{\bco}{\begin{corollary}}
\newcommand{\eco}{\end{corollary}}
\newcommand{\bpro}{\begin{proposition}}
\newcommand{\epro}{\end{proposition}}

\newtheorem{assumption}{Assumptions}[theorem]
\newcommand{\bas}{\begin{assumption}}
\newcommand{\eas}{\end{assumption}}

\theoremstyle{definition}

\newtheorem{example}[theorem]{Example}
\newtheorem{remark}[theorem]{Remark}
\newtheorem{remarks}[theorem]{Remarks}
\newtheorem*{proofNEW}{Proof}
\newcommand{\bpr}{\begin{proofNEW}}
\newcommand{\epr}{\qed
\medskip
\end{proofNEW}}

\theoremstyle{plainNEW}

\newcommand{\bth}{\begin{theorem}}
\renewcommand{\eth}{\end{theorem}}
\newtheorem{proposition}[theorem]{Proposition}

\pagenumbering{gobble}
\title{Stability of the  Solitary Manifold\\ 
of the Perturbed Sine-Gordon Equation}

\date{}

\author{{\sc Timur Mashkin}\\[2ex]
         Mathematisches Institut, Universit\"at K\"oln, \\
         Weyertal 86-90, D\,-\,50931 K\"oln, Germany \\
         e-mail: tmashkin@math.uni-koeln.de}

\maketitle

\begin{abstract}\noindent
We study the perturbed sine-Gordon equation 
$\theta_{tt}-\theta_{xx}+\sin \theta=  F(\eps,x)$,
where $F$ is of differentiability class $C^n$ in $\eps$ and the first $k$ derivatives vanish at $0$, i.e., $\deps^l F(0,\cdot)=0$ for $0\le l\le k $. 
We construct implicitly a virtual solitary manifold by deformation of the classical solitary manifold in $n$ iteration steps.
Our main result establishes that the initial value problem with an appropriate initial state $\eps^n$-close to the virtual solitary manifold has a unique solution which follows up to time $1/(\ti C{\eps^{\fr{k+1}2}})$ and errors of order $\eps^n$ a trajectory on the virtual solitary manifold. The trajectory on the virtual solitary manifold is described by two parameters which satisfy a system of ODEs. 
In contrast to previous works our stability result yields arbitrarily high accuracy as long as the perturbation $F$ is sufficiently often differentiable.  
%
\end{abstract}

\pagenumbering{arabic}

\section{Introduction}

\noindent 
The perturbed sine-Gordon equation 
\beq\la{SGE}
\t_{tt}-\t_{xx}+\sin\t=F(\eps,x),~~~~t,x\in\R,~~~~\eps\ll 1,
\eeq
is a Hamiltonian evolution equation with Hamiltonian given by
\beq
H^\eps(\t,\p)=\fr 1 2\int \p^2+\t_x^2+2(1-\cos\t)-2F(\eps,x)\t \,dx\,\nn
\eeq
and the symplectic form given by
\be\la{symplecticform introduction}
\Omega\l(\tpP,\tp \r) = \l\langle \tpP,\J\tp \r\rangle_{L^2(\R)\oplus L^2(\R)}=\int_\R \p'(x)\t(x)-\t'(x)\p(x)\,dx,
\ee
where 
$$\J=\l(\begin{matrix}
0&-1\\
1&0\\
\end{matrix}\r).$$ 
In first order formulation \re{SGE} can be written as a system:
\be\la{SGE1 first order introduction}
\partial_t\bma
\t\\
\p
\ema=\l(\begin{matrix}
\p\\
\t_{xx}-\sin\t+F(\eps,x)\\
\end{matrix}\r).
\ee
The unperturbed sine-Gordon equation ($F(\eps,x)=0$),
admits soliton solutions 
$
\bma 
\t_0(\xi(t),\u(t),x)\\
\p_0(\xi(t),\u(t),x)
\ema
$,
where
$$
\dot\xi=\u\,, ~~ \dot \u=0\,,~~~~(\xi(0),\u(0))=(a,v)\in\R\times(-1,1).
$$
Here the functions $(\t_0,\p_0)$ are defined by
\be\la{solitonsolution}
{}&\bma
\t_0(\xi,\u,x)\\
\p_0(\xi,\u,x)
\ema:=\bma
\t_K(\g(u)(x-\xi))\\
-\u\g(u)\t_K'(\g(u)(x-\xi))\\
\ema\,,~\u\in(-1,1),~~\xi,x\in\R,
\ee
where 
$$
\g(u)=\frac 1 {\sqrt{1-u^2}},~~~~\t_K(x) =4\arctan(e^x),
$$
and $\t_K$ satisfies $\t_K''(x)=\sin\t_K(x)$ with boundary conditions $\t_K(x) \to \bma 2\pi\\0  \ema$ as $x\to \pm \infty$.
The states $\l(\begin{matrix}
\t_0(a,v,\cdot)\\
\p_0(a,v,\cdot)\\
\end{matrix}\r)$ form the classical two-dimensional solitary manifold  
$$
{\cal S}_0:=\l\{ \l(\begin{matrix}
\t_0(a,v,\cdot)\\
\p_0(a,v,\cdot)\\
\end{matrix}\r)~:~v\in(-1,1),~a\in\R
\r\}.
$$
In the present paper, we assume that the perturbation term $F$ in \re{SGE} is of differentiability class $C^n$ in $\eps$ and that its first $k$ derivatives vanish at 0, i.e., 
\be\la{intro assumption on F}
\deps^l F(0,\cdot)=0~~ \text{for} ~~0\le l\le k.
\ee
We construct 
a virtual solitary manifold ${\cal S}_n^\eps$, which is adjusted to the perturbation $F$. 
The construction can be thought of as a successive distortion of the classical solitary manifold ${\cal S}_0$. It is based on an iteration scheme composed
of n steps, where in each iteration step a specific PDE will be solved implicitly.
In the last iteration step we obtain an implicit solution  $(\t_n^\eps(\xi,u,x),
\p_n^\eps(\xi,u,x))$ which defines
the virtual solitary manifold  
\be\la{virtMF}
{\cal S}_n^\eps:=\l\{ \bma
\t_n^\eps(a,v,\cdot)\\
\p_n^\eps(a,v,\cdot)
\ema~:~v\in(-u_*,u_*),~a\in\R
\r\}, ~~~~u_*\in (0,1].
\ee
   %
%
We consider for $\xi_s\in\R$ and $\eps\ll 1$ the Cauchy problem
\be\la{Cauchy_intro}
\partial_t \bma
\t \\
\p 
\ema =\bma
\p \\
\dx^2\t -\sin\t +\F 
\ema,~~
\bma
\t(0,x)\\
\p(0,x)
\ema =
\bma
\t^\eps_n(\xi_s,\u_s,x)\\
\p^\eps_n(\xi_s,\u_s,x)
\ema+
\bma
\v(0,x)\\
\w(0,x) 
\ema,
\ee
with initial data $\eps^n$-close to the virtual solitary manifold ${\cal S}_n^\eps$, i.e., 
$$
\nhone{v(0,\cdot)}^2+\nw{w(0,\cdot)}^2\le \eps^{2n}.
$$ 
Further, we suppose that 
$(\v(0,\cdot),\w(0,\cdot))$ is symplectic orthogonal to the tangent space of ${\cal S}_n^\eps$ at  
$(\t^\eps_n(\xi_s,\u_s,\cdot),
\p^\eps_n(\xi_s,\u_s,\cdot))$
and that the smallness assumption
$$
|u_s|\le \tilde C\eps^{\fr{k+1}2}
$$
is satisfied.

Our main theorem shows that, under the mentioned assumptions on the initial data, the Cauchy problem \re{Cauchy_intro} has a unique solution $(\t,\p)$ for times 
\be\la{intro timescale}
0 \le t\le \fr 1 { \tilde C \eps ^{\fr{k+1}2}},
\ee 
which may be written in the form
\be
\bma
\t(t,x)\\ 
\p(t,x)
\ema=\bma 
\t_n^\eps(\bar\xi(t),\bar\u(t),x)\\ 
\p_n^\eps(\bar\xi(t),\bar\u(t),x)
\ema+ 
\bma
\v(t,x)\\
\w(t,x)
\ema.\nn
\ee
Furthermore, the solution remains $\eps^n$-close to the manifold ${\cal S}_n^\eps$, i.e.,
\be\la{bound on vw introduction}
\nhone{v(t,\cdot)}^2+\nltwo{w(t,\cdot)}^2\le \tilde C \eps^{2n}
\ee
and the parameters $(\bar\xi(t),\bar\u(t))$ satisfy the ODEs
\be\la{ODE introduction}
 \d{\bar\xi}( t) =  \bar\u(t)  \,,~~~~
 \d{\bar\u}( t) = \lambda_{n}^\eps\l(\bar\xi(t), \bar\u(t)\r), 
\ee
with initial data 
$
\bar\xi(0)=\xi_s,~\bar\u(0)=\u_s
$, where $\lambda_{n}^\eps$ is defined implicitly. The time scale \re{intro timescale}
is nontrivial and the 
parameters $\bar\xi,\bar\u$ describe
a fixed nontrivial perturbation of the uniform linear motion
as $\eps \to 0$ if the perturbation $F$  satisfies condition \re{intro condition on nontrivial dynamic} mentioned below.
\\\\
This result yields a fairly accurate description of the solution $(\t,\p)$ to the Cauchy problem \re{Cauchy_intro}, since we are able to control the dynamics on the virtual solitary manifold ${\cal S}_n^\eps$  by the ODEs \re{ODE introduction} and the dynamics of the transversal component $(v(t,\cdot),w(t,\cdot))$ by the upper bound on its norm \re{bound on vw introduction}. 
\\
The higher the differentiability class $C^n$
of the perturbation $F$ the higher is the accuracy of our stability statement. 
The time scale of the result is the larger
the more first derivatives of $F$ vanish at 0.
A precise statement is found in \cref{ch: Main Results virtual}.
\\
\\
Let us mention related works and give some background to our paper. 
Orbital stability of soliton solutions under perturbations of the initial data has been proven for the (unperturbed) sine-Gordon equation (see \cite{MR678151}, \cite[Section 4]{Stuart3}).
In \cite{Stuart2} D. M. Stuart investigated the perturbed sine-Gordon equation,
\be
{}&\t_{TT}-\t_{XX}+\sin\t+\eps g=0,\nn
\ee
where the perturbation $g=g(\t)$ is a smooth function such that $g_0(Z)=g(\t_K(Z))\in L^2(dZ)$ and $\eps\ll 1$. He proved  that there exists $T_* = {\cal O} \l( \fr 1 \eps \r) $ such that the corresponding initial value problem with initial data
\be
{}&\t(0,X)=\t_K(Z(0))+\eps \ti\t(0,X),~~~~
\t_T(0,X)=\fr {-u(0)}{\sqrt{1-u(0)^2}}\t_K'(Z(0))+\eps \ti\t_T(0,X),\nn\\
{}& (\ti\t(0,X),\ti\t_T(0,X))\in H^1\oplus L^2,\nn
\ee
has a unique solution of the form 
\be
\t(T,X)=\t_K(Z)+\eps\ti\t(T,X),~~~Z=\fr{X-\int^T u - C(T)}{\sqrt{1-u^2}},\nn
\ee
where $\ti\t\in C([0,T_*],H^1),~\t_T\in C([0,T_*],L^2)$ and
\be
{}&C(T)=C_0(\eps T)+\eps \ti C,~~
{}&u(T)=u_0(\eps T)+\eps \ti\u(T) \l( \Rightarrow~ p= \fr {u} {\sqrt{1-u^2}}= p_0(\eps T)+\eps \ti p(T)  \r).\nn
\ee
Here $\ti p,~\ti u,\ti C, \fr {d\ti u}{d T}, \fr {d\ti C}{d T}, |\ti\t |_{H^1(\R)}$ are bounded independent of $\eps$ and $u_0,C_0$ are solutions of certain explicitly given modulation equations.
The proof is based on an orthogonal
decomposition of the solution into an oscillatory part and a one-dimensional
"zero-mode" term.
 
The sine-Gordon equation arises in various physical phenomena such as dynamics of long Josephson junctions \cite{0953-8984-7-2-013,RevModPhys.61.763}, dislocations in crystals 
\cite{FrenKont}, waves in ferromagnetic materials \cite{0022-3719-11-1-007}, etc. In \cite{Skyrme237} T. H. R. Skyrme proposed the equation to model elementary particles. Dynamics of solitons under constant electric field were examined numerically in \cite{doi:10.1143/JPSJ.46.1594}.
In the present paper, we investigate virtual solitons in the presence of a time independent electric field $F(\eps,x)$, which is a physically relevant problem.

There are also many other long (but finite)-time results for different equations with external potentials such as \cite{MR2094474,MR2232367,MR2342704,MR2855072}.
For instance, in \cite{HoZwSolitonint} J. Holmer and M. Zworski considered the Gross-Pitaevskii equation 
\be 
\left\{
\begin{aligned}
&
i\partial_t u + \tfrac{1}{2}\partial_x^2 u - V ( x ) u +u|u|^2 = 0,
\\
&
u(x,0) = e^{ i v_0 x } \sech ( x - a_0 ),  
\end{aligned}\nn
\right.
\ee
with a slowly varying smooth potential $V(x)=W(\eps x)$, where $W\in C^3(\R,\R)$. 
They proved that up to time $\fr{\log(1/\eps)}\eps$ and errors of size $\eps^2$ in $H^1$, the solution is a soliton evolving according to the classical dynamics of a natural effective Hamiltonian. 
 
A common method for investigating the interaction of solitons with external potentials is 
to decompose 
the dynamics 
in a neighbourhood of the classical solitary manifold and to
apply Lyapunov-type arguments afterwards.
Beside other main ingredients this approach has been chosen in \cite{MR2094474,MR2232367,MR2342704,HoZwSolitonint,MR2855072}.  

In the present paper we extend this method for the perturbed sine-Gordon equation by
introducing, the adjusted to the perturbation $F $, virtual solitary manifold.
We decompose the dynamics in a part on the virtual (rather than classical) solitary manifold and a transversal part, before we proceed with a customized Lyapunov method. 
Utilizing the virtual solitary manifold is crucial for the high accuracy in our stability result.

In the broadest sense, a similar technique has been used in \cite[Section 4, Section 5]{HolmerLin} for the NLS equation.
In the present paper,
in the iteration scheme for the construction of the virtual solitary manifold, each iteration corresponds to a correction of the classical solitary manifold. 
These corrections are represented by approximate solutions of a specific PDE \re{intro spec PDE} mentioned below, 
where the accuracy of the approximate solutions increases with each iteration step by order $1$ in $\eps$.
So in \cite[Section 4, Section 5]{HolmerLin},
the solitary manifold has been corrected   once, which corresponds to the correction executed in the first iteration in the construction of the virtual solitary manifold in the present paper.
In our approach, the construction of the virtual solitary manifold relies on the implicit function theorem.
However, the correction in \cite[Section 4, Section 5]{HolmerLin} was done in the form
of a direct asymptotic expansion 
and not by employing the implicit function theorem.
In that sense our paper presents a new point of view. 
 
A further remarkable point is that unlike \cite[Section 4, Section 5]{HolmerLin} we are able to carry out arbitrarily many corrections in the case of the sine-Gordon equation 
as long as the perturbation $F$ is sufficiently often differentiable in $\eps$.
The accuracy of our stability result is the higher the more corrections are possible, wheras the number of possible corrections is determined by the differentiability class of the perturbation $F$.
\\ 
We abstained from considering perturbations of type $\eps W(\eps x)$ for the following reason. In our approach we do need the assumption that the perturbation $F$ is  differentiable with respect to $\eps$, but
there does not exist a function $W\not= 0,W\in L^2(\R)$
such that the mapping 
$
\eps\mapsto\eps W(\eps\cdot)\,
$
is differentiable in $L^2(\R)$.
%

Further results on long time soliton asymptotics and orbital stability for different equations can be found in
\cite{MR820338,MR0428914,MR0386438,MR2920823,MR1071238,MR1221351,ImaykinKomechVainberg,MR3630087,MR3461359}.
This paper is based on \cite[Part IV]{Mashkin}, where many of the computations are presented in greater detail.

Now let us comment on our techniques. The virtual solitary manifold is constructed by the following iteration scheme. 
Let firstly $\eps\mapsto \tilde F(\eps)$ be a general function of differentiability class $C^n$ 
mapping into a specific 
Sobolev space
such that $\tilde F(\eps)$ depends on $(\xi,x)$
and 
$\tilde F(0)=0$.
$\tilde F$ will be specified later.  
The function $(\t_0,\p_0)$, given by \re{solitonsolution}, is a solution of 
\be
{}&\BR{\u\dxi\l(\begin{matrix}
\t\\
\p\\
\end{matrix}\r)
-\l(\begin{matrix}
\p\\
\dx^2\t-\sin\t\\
\end{matrix}\r)
}{\large$=:{\cal G}_0(\t,\p)$}=0,\nn
\ee
which is the equation characterizing the classical solitons. 
In the first iteration step we modify 
${\cal G}_0(\t,\p)=0$ 
by introducing an
additional unknown variable $\lambda$ and 
adding some terms involving $(\t_0,\p_0)$ and  $\tilde F $. 
The modified equation is of the form
\be\la{successive eq G1}
{}&\BR{\u\dxi\bma
\t\\
\p\\
\ema
-\l(\begin{matrix}
\p\\
\t\xx-\sin\t+\ti F (\eps)\\
\end{matrix}\r)
+\lambda \du\bma
\t_0\\
\p_0\\
\ema
}{\large$=:\Go(\t,\p,\lambda)$}=0\,.
\ee
Here and in the following iterations the functions $\t,\p$ depend on $(\xi,\u,x)$ and
$\lambda$ depends on $(\xi,\u)$. 
We solve ${\cal G}_1^\eps(\t,\p,\lambda)=0 $ implicitly for $(\t,\p,\lambda)$ in terms of $\eps$ and denote the solution by $(\t_1^\eps,\p_1^\eps,\lambda_{1}^\eps)$. 
%
In the next iteration step we modify ${\cal G}_1^\eps(\t,\p,\lambda)=0 $
by adding some terms involving  $(\t_1^\eps,\p_1^\eps)$ and solve the modified equation of the form 
\be\la{successive eq G2}
{}&\BR{\u\dxi\bma
\t\\
\p\\
\ema
-\l(\begin{matrix}
\p\\
\t\xx-\sin\t+\ti F(\eps)\\
\end{matrix}\r)
+\lambda\du\bma
\t_1^0+\deps\t_1^0\eps\\
\p_1^0+\deps\p_1^0\eps\\
\ema 
}{\large$=:\Gt(\t,\p,\lambda)$}=0\,
\ee
implicitly for $(\t,\p,\lambda)$ in terms of $\eps$. 
Due to the assumption that $\eps\mapsto \tilde F(\eps)$
is of differentiability class $C^n$, 
it is possible to iterate this 
procedure until we obtain in the $n$th step
the equation
\be\la{successive eq Gn}
{}&\BR{\u\dxi\bma
\t\\
\p\\
\ema
-\l(\begin{matrix}
\p\\
\t\xx-\sin\t+\ti F(\eps)\\
\end{matrix}\r)
+\lambda\du\bma
\sum_{i=0}^{n-1} \fr{\deps^i\t_{n-1}^0}{i!}\eps^i\\
\sum_{i=0}^{n-1} \fr{\deps^i\p_{n-1}^0}{i!}\eps^i\\
\ema 
}{\large$=:\Gn(\t,\p,\lambda)$}=0.
\ee
We solve $\Gn(\t,\p,\lambda)=0$ implicitly for $(\t,\p,\lambda)$ in terms of $\eps$ and 
denote the solution by $(\t_n^\eps,\p_n^\eps,\lambda_{n}^\eps)$.  
The
existence of the implicit solutions $\eps\mapsto(\t_j^\eps,\p_j^\eps,\lambda_{j}^\eps)$ for $1\le j\le n$  
is ensured by the implicit function theorem.
In the actual proof,
we consider for functional analytic reasons the translated maps and solve the equations 
\be\la{intro def tilde G}
\ti {\cal G}_j^\eps(\hatt,\hatp,\lambda):={\cal G}_j^\eps(\t_0+\hatt,\p_0+\hatp,\lambda)=0
\ee
for $ (\hatt,\hatp,\lambda)$ in terms of $\eps$. This is caused, among others, by the fact that $\t_0(\xi,\u,x)  \not\rightarrow 0$ as $|x|\to \infty$ for fixed $\xi$ and $\u$. 
We denote the solutions to the equations $\ti {\cal G}_j^\eps(\hatt,\hatp,\lambda)=0$ by $( \hatt_{j}^\eps , \hatp_{j}^\eps ,\lambda_{j}^\eps)$, where $(\t_{j}^\eps ,\p_{j}^\eps ,\lambda_{j}^\eps)
 =(\t_0+\hatt_{j}^\eps ,\p_0+\hatp_{j}^\eps ,\lambda_{j}^\eps)$. 
The application of the implicit function theorem
relies on the fact that 
$(0,0,0,0)$ 
solves all equations in a particular point, i.e., 
$\ti {\cal G}_j^0(0,0,0)=0$. 
As a consequence of the construction, the 
solution obtained in the $j$th iteration
$\eps\mapsto(\t_j^\eps,\p_j^\eps,\lambda_{j}^\eps)$ solves the equation 
\be\la{intro spec PDE}
{\u\dxi\bma
\t\\
\p\\
\ema
-\l(\begin{matrix}
\p\\
\t\xx-\sin\t+\ti F (\eps)\\
\end{matrix}\r)
+\lambda \du\bma
\t\\
\p\\
\ema  
}
=0\,
\ee
up to errors of order $\eps^{j+1}$ for $1\le j\le n$.
 
In order to define the virtual solitary manifold we apply this iteration scheme on a specific $\ti F$,  which is
a truncated version of the perturbation term $F$ in \re{Cauchy_intro}, given by
\be\la{intro assumption Chi and F}
\begin{cases}
 \ti F(\eps,\xi,x):=\F \chi(\xi),\\ 
 \text{where } \chi\in C^{\infty}(\R),~\chi (\xi)=1 \text{  for } |\xi|\le |\xi_s|+3 \text{ and } \chi (\xi)=0 \text{ for } |\xi|\ge |\xi_s|+4.
\end{cases}
\ee
From now on we denote by $(\t_n^\eps,\p_n^\eps,\lambda_n^\eps)$ the solution obtained in the $n$th iteration by application of the iteration scheme on the specific $\ti F$ given by \re{intro assumption Chi and F}.  
The first two components of 
$(\t_n^\eps,\p_n^\eps,\lambda_n^\eps)$ define the virtual solitary manifold \re{virtMF}.  
A further consequence of the construction and assumption \re{intro assumption on F} is that the functions
$(\t^\eps_n(\xi(t),\u(t),x), 
\p^\eps_n(\xi(t),\u(t),x)) $
solve the perturbed sine-Gordon equation \re{SGE1 first order introduction}
up to an error of order $\eps^{n+k+1}$ as long as $(\xi(t),\u(t))$ satisfy the ODE system 
$\dot\xi(t)=u(t)$, $\dot \u(t)=\lambda^\eps_n(\xi(t),\u(t))$. 
We call these approximate solutions virtual solitons. In the further proof they play a role which is comparable to that of classical solitons, for instance, in the proof of orbital stability ($F=0$) of classical solitons (see \cite[Section 4]{Stuart3}). 

The idea of deforming the classical solitary manifold and
utilizing thereby implicitly defined functions appears in \cite{Stuart3} with the purpose of rewriting the Hamiltonian in a neighbourhood of the manifold of virtual solitons (see \cite[Section 3]{Stuart3}). 
The virtual solitons in our paper
and the corresponding manifold \re{virtMF} are defined by 
equations and an iteration scheme that were not considered in \cite{Stuart3}.

The existence of a local solution of the Cauchy problem \re{Cauchy_intro} follows from the contraction mapping theorem. 
In the following approach we
derive some bounds which 
imply that the local solution is continuable and that estimate \re{bound on vw introduction} is satisfied on the relevant time scale. 
We decompose the solution of \re{Cauchy_intro}
into a point on the virtual solitary manifold ${\cal S}_n^\eps$ and a transversal component, 	
i.e.,
\be\la{intro_decomposition}
\bma
\t(t,x)\\ 
\p(t,x)
\ema=\bma 
\t_n^\eps( \xi(t), \u(t),x)\\ 
\p_n^\eps( \xi(t), \u(t),x)
\ema+ 
\bma
\v(t,x)\\
\w(t,x)
\ema,
\ee
where the parameters $(\xi(t),\u(t))$ are chosen in such a way that the transversal component $(v(t,\cdot),w(t,\cdot))$ is symplectic orthogonal to the tangent space of ${\cal S}_n^\eps$ at the corresponding point.
This symplectic decomposition is possible in a small uniform distance to the virtual solitary manifold due to the implicit function theorem.
The energy
\be
{}& H(\t,\p)=\fr 1 2\int \p^2+\t_x^2+2(1-\cos\t) \,dx\nn
\ee
and the momentum 
\be
{}&\Pi (\t ,\p )
 = \int\p\t_x\,dx \nn
\ee
are conserved quantities of the unperturbed sine-Gordon equation. We make use of this fact and achieve control over the transversal component of the solution $(v,w)$ by utilizing an almost conserved  
Lyapunov function, given by
\be 
 L^\eps =\int \fr{\w^2 } 2 +\fr{ \dx\v^2} 2 +\fr{\cos(\t_n^\eps(\xi,\u,\cdot)) \v^2 }2+\u\w \dx\v \,dx,\nn
\ee
where $(v,w)$ and $(\xi,\u)$ are such as in \re{intro_decomposition}. $ L^\eps$ is the quadratic part of
$$
H(\t_n^\eps +v,\p_n^\eps+w)+u\Pi(\t_n^\eps+v,\p_n^\eps+w) - \Big (H(\t_n^\eps,\p_n^\eps)+u\Pi(\t_n^\eps,\p_n^\eps) \Big).
$$
The Lyapunov function is bounded from below in terms of $\nhone{v(t,\cdot)}^2+\nltwo{w(t,\cdot)}^2$, 
which is a consequence of symplectic orthogonality in decomposition \re{intro_decomposition} and of 
spectral properties of the operator
$-\dZ^2+\cos\t_K(Z)$. 
The parameters $(\xi,\u)$ 
satisfy ODEs \re{ODE introduction} up to errors of order $\eps^{n+k+1}$, which goes 
back (among others) to the
construction of the virtual solitary manifold, especially to the fact that $(\t_n^\eps,\p_n^\eps,\lambda_n^\eps)$ solves \re{successive eq Gn} with $\tilde F$ given by \re{intro assumption Chi and F}. 
This property of 
$(\xi,\u)$ and once again the mentioned fact about $(\t_n^\eps,\p_n^\eps,\lambda_n^\eps)$ allow us to 
control the Lyapunov function from above.
Therefore we are able to estimate the norm of the transversal component $(v,w)$ and obtain ultimately bound \re{bound on vw introduction}.
Using Gronwall's lemma we pass from the approximate equations for the parameters $(\xi,\u)$ to the exact ODEs \re{ODE introduction}.

Finally let us explain
under which conditions
the result provides a nontrivial dynamics on the virtual solitary manifold 
as $\eps \to 0$. 
The linearization of 
$(\hatt,\hatp,\lambda) \mapsto 
\ti {\cal G}_n^\eps(\hatt,\hatp,\lambda)
$ 
carried out at 
$(\hatt,\hatp,\lambda)=(0,0,0)$, $\eps=0$ 
is invertible
and we denote the linearization by
 $$
{\frak M}_2^1:  
( \t,
\p,
\lambda) \mapsto
{\frak M}_2^1
( \t,
\p,
\lambda ).
$$
Thus there exist functions $( \bar\t,
\bar \p,
\bar \lambda )$ such that the $(k+1)$th derivative of $\ti F$ with respect to $\eps$, evaluated at $\eps=0$, can be written in the form
\be \la{intro condition on tiF}
\bma
0\\
 \deps^{k+1}  \ti F(0)
\ema
=
{\frak M}_2^1
( \bar\t,
\bar \p,
\bar \lambda ),~~~
%
\text{  ${\frak M}_2^1$ given by \cref{le invertibilityMxiCtwo alpha} ($n=2,\a=1$)}.
\ee
Here the functions $\bar\t,\bar\p$ depend on $(\xi,\u,x)$ and $\bar\lambda$ depends on $(\xi,\u)$.
The ODEs \re{ODE introduction} can be rescaled in time by introducing $s=\eps^{\beta(k )} t$ with $\beta(k )=\fr{k+1} 2$,
$
\hxi(s)= \bar\xi(s/\eps^{ \beta(k ) })
$,
and
$
\hu(s)= \fr 1 {\eps^{\beta(k ) }} {\bar u(s/\eps^{ \beta(k )})}
$ 
such that the 
corresponding transformed ODEs have the form 
\be
\ds \hxi(s) =   \hu(s) , ~~~~
\ds \hu(s) =  \fr 1 {\eps^{2\beta(k )}} \lambda_{ n}^\eps(\hxi(s), \eps^{\beta(k )}\hu(s)).\nn
\ee
As $\eps \to 0$, the transformed ODEs converge to ODEs that describe a fixed nontrivial perturbation of the uniform linear motion if 
the next condition is satisfied:
\be
\la{intro condition on nontrivial dynamic}
\begin{cases}
{}&\text{There exists $\chi$ satisfying \re{intro assumption Chi and F} such that for $\ti F$ given by \re{intro assumption Chi and F} the following}\\ {}& \text{holds: The function $\bar \lambda$ in 
representation \re{intro condition on tiF} fulfills }\bar \lambda(\cdot,0)\not= 0.  
\end{cases}
\ee 
This 
is for the following reason.
Due to \re{intro assumption on F} $\deps^l  \ti F(0)=0$ for $1\le l\le k$ and differentiation of ${\cal G}_n^\eps(\t_n^\eps,\p_n^\eps,\lambda_n^\eps)=0$ 
with respect to $\eps$ yields (see proof of \cref{thITrelations}):   
\be 
\bma
0\\
 \deps^l  \ti F(0)
\ema
=
{\frak M}_2^1
( \deps^l \t_n^0,
\deps^l \p_n^0,
\deps^l\lambda_n^0 ),~~~~~1\le l\le k+1.
\ee
Using invertibility of ${\frak M}_2^1$, condition \re{intro condition on nontrivial dynamic} and the fact that $\lambda_n^0=0$ 
it follows that $0\not=\lambda_n^\eps= {\cal O}(\eps^{k+1})$, which implies the claim. 

The paper is organized as follows. In \cref{ch: Main Results virtual}, we formulate the main result. In \cref{Virtual Solitary Manifold}, we construct the virtual solitary manifold. 
We prove in \cref{ch:Symplectic Orthogonal Decomposition} that in a uniform distance to the virtual solitary manifold a decomposition into symplectically orthogonal components is possible. The existence of a local solution $(\t,\p)$ with initial state close to the virtual solitary manifold is established in \cref{local solution F virtual}. In \cref{ch Modulation Equations}, we derive modulation equations for the parameters that describe the position on the manifold. We introduce a Lyapunov function and compute its time derivative in \cref{ch:Lyapunov functional}. A lower bound on the Lyapunov function is proved in \cref{ch: Lower bound}. In \cref{aprioriestimate F virtual manifold}, we prove our main result, \cref{maintheorem specialF kderivativesvanish virtual manifold}.
Some preliminary decompositions are showed in Appendix A. These decompositions are used in Appendix B, where we prove that the  linearizations considered in \cref{Virtual Solitary Manifold} are invertible.  
\subsection*{Notation and Conventions}
For a Hilbert space $H$ its inner product is denoted by $\langle\cdot,\cdot\rangle_H $, the orthogonal complement of a closed subspace $M$ 
in $H$ by $M^{\perp, H}$, the orthogonal projection on $M$ by $(\cdot)_M$ and the span of $v_1,\ldots, v_p\in H$ by $\langle v_1,\ldots, v_p\rangle$ .
For functions $\lambda$ depending on $(\xi,\u)$ and functions $\t$ depending on $(\xi,\u,x)$ the notation $\lambda(\xi,\u)=\lambda(\u)(\xi)$, $\t(\xi,\u,x)=\t(\u)(\xi,x)$ is used. $\g$ without an argument denotes always $\g(u)$. Occasionally we drop the dependence of functions on certain variables. We also denote occasionally by $\Vert \cdot \Vert $ the norm of an operator and drop the spaces in the notation. We write ${L_x^{2}(\R)},{H_{\xi,x}^{k}(\R^2)}$ and so on for the Lebesgue and Sobolev spaces when we wish to emphasize the variables of integration.   
\section{Main Result}\la{ch: Main Results virtual}
To formulate our result precisely, we need some definitions.
\bde\la{def:PartfourMainResult}
Let $\a,n\in\N$ and $u_*>0$. Let us denote by $I(\u_*):=[-\u_*,\u_*]$.
\begin{itemize}
\item [(a)] $H^{k,\alpha}(\R) $ denotes the weighted Sobolev space of functions with finite norm
$$|\t|_{H^{k,\alpha}(\R)}= |(1+|x|^2) ^\fr \alpha 2\t(x)|_{H_x^{k}(\R)}.$$
\item [(b)] $H^{k,\alpha}(\R^2) $ denotes the weighted Sobolev space of functions with finite norm
$$|\t|_{H^{k,\alpha}(\R^2)}= |(1+|\xi|^2+|x|^2) ^\fr \alpha 2\t(\xi,x)|_{H_{\xi,x}^{k}(\R^2)}.$$
\item [(c)] $ \ubar{ Y}^\a$ is the space $H^{3,\a}(\R^2) \oplus H^{2,\a}(\R^2) \oplus H^{2,\a}(\R)$
with the finite norm
$$
|y|_{\ubar{ Y}^\a} = |\t|_{H^{3,\a}(\R^2)}+ |\p|_{H^{2,\a}(\R^2)}+|\lambda |_{H^{2,\a}(\R)}.
$$ 
\item [(d)] $Y_n^\a(\u_*)$ is the space\\
$\\\ba
{}&\bigg\{ y=(\t,\p,\lambda) \in C^n( I(\u_*), \ubar{ Y}^\a) :  \Vert y \Vert_{Y_n^\a(\u_*)} <\infty;~\forall~ \u\in I(\u_*),~\forall~\mu\in H^{2,\a}(\R):\\
{}& \Ltwortwoaxix{ \bma \t(\u)(\xi,x)\\
\p(\u)(\xi,x) \ema}{\mu(\xi)\bma \t_K'(\Z)\\
-\u\g\t_K''(\Z)\ema} =0 \bigg\}\,
\ea\\$\\
with the finite norm
$$
\Vert y \Vert_{Y_n^\a(\u_*)} =\sup_{\u\in I(\u_*)}  \l( \sum_{i=0}^n |\du^i y(u)|_{\ubar{ Y}^\a}\r).
$$
\item[(e)] For $l\in {\mathbb N}$ and $0<U<u_*$ we introduce the parameter area
\be
\Sigma(l,U,u_*):=\Big\{(\xi,\u)\in \R\times (-1,1): 
u\in  (-U- V(l,U,u_*),U+ V(l,U,u_*) ) \Big\},\nn
\ee
where $V(l,U,u_*):= \fr{u_*-U} l$. 
\end{itemize}
\ede

\noindent The weighted Sobolev spaces in \cref{def:PartfourMainResult} (a), (b) are defined as in \cite{Kopylova}. We are now ready to state our main result.
\bth\la{maintheorem specialF kderivativesvanish virtual manifold}
Let $n,k\in{\mathbb N}$,  
$n\ge 1$, $k+1 \le n$.
Assume that $\xi_s\in \R $, $F\in C^{n}((-1,1),H^{1,1}(\R))$ and $\deps^l F(0,\cdot)=0$ for  $0\le l\le k $.
Then there exist $\eps_0, \u_*, \tilde C>0$ and a map 
$$
(-\eps_0,\eps_0) \to Y_{2}^1(\u_*),~
\eps \mapsto (\hat\t_n^\eps,\hat\p_n^\eps,\lambda_{n}^\eps)\la{map in main theorem}
$$
of class $C^n$ such that the following holds. 
Let $\eps\in(0,\eps_0)$ and $0<U<u_*$. Consider the Cauchy problem
\be\la{SGE1}
\partial_t \bma
\t \\
\p 
\ema
=\l(\begin{matrix}
\p \\
\dx^2\t -\sin\t +\F\\
\end{matrix}\r),
\bma
\t(0,x)\\
\p(0,x)
\ema=\bma
\t^\eps_n(\xi_s,\u_s,x)\\
\p^\eps_n(\xi_s,\u_s,x)
\ema 
+ \bma\v(0,x)\\
  \w(0,x)
\ema,
\ee
where $(\t_n^\eps,\p_n^\eps,\lambda_{n}^\eps)=(\t_0+\hatt_n^\eps ,\p_0+\hatp_n^\eps ,\lambda_{n}^\eps)$ with $(\t_0,\p_0)$ given by \re{solitonsolution} 
such that the following assumptions are satisfied:
\begin{itemize} 
\item[(a)] $|u_s|\le \tilde C\eps^{\fr{k+1}2}$;
\item[(b)] ${\cal N}^\eps(\t(0,\cdot),\p(0,\cdot),\xi_s,\u_s)=0$,\\ 
where 
$
{\cal N}^\eps=({\cal N}^\eps_1,{\cal N}^\eps_2) : L^\infty(\R) \times L^2(\R)\times \Sigma(2,U,u_*)\to\R^2 
$ is given by
\be\la{orthogonalitycondIntro}
{\cal N}^\eps(\t,\p,\xi,\u)
:=
\bma
\OM{\bma
\dxi\t_n^\eps(\xi,\u,\cdot)\\
\dxi\p_n^\eps(\xi,\u,\cdot)
\ema}{\bma
\t(\cdot)-\t_n^\eps(\xi,\u,\cdot)\\
\p(\cdot)-\p_n^\eps(\xi,\u,\cdot)\\
\ema}\\
\OM{\bma
\du\t_n^\eps(\xi,\u,\cdot)\\
\du\p_n^\eps(\xi,\u,\cdot)\\
\ema}{\bma
\t(\cdot)-\t_n^\eps(\xi,\u,\cdot)\\
\p(\cdot)-\p_n^\eps(\xi,\u,\cdot)\\
\ema}
\ema\,
\ee
and the symplectic form $\Omega$ is given by \re{symplecticform introduction};
\item[(c)] $\nhone{v(0,\cdot)}^2+\nltwo{w(0,\cdot)}^2\le \eps^{2n}$, where $(\v(0,\cdot),\w(0,\cdot))$ is given by 
\re{SGE1}.
\end{itemize}
Then the Cauchy problem defined by \re{SGE1}
has a unique solution on the time interval
\be
0\le t \le T,  ~\text{where}~ T=T(\eps,k)=\fr {1} {\tilde C\eps^{\beta(k )}}, ~~\beta(k )=\fr{k+1 }2.\nn
\ee
The solution may be written in the form
\be
\bma
\t(t,x)\\ 
\p(t,x)
\ema=\bma 
\t_n^\eps(\bar\xi(t),\bar\u(t),x)\\ 
\p_n^\eps(\bar\xi(t),\bar\u(t),x)
\ema+ 
\bma
\v(t,x)\\
\w(t,x)
\ema,\nn
\ee
where $\v, \w,$ 
have regularity
$
(v(t), w(t)) \in C^1([0,T ] , H^1(\R) \oplus L^2(\R))
$
and $\bar\xi,\bar\u$ solve the ODEs
\be\la{exactODE virtual1}
 {}&\d{\bar\xi}( t) =  \bar\u(t)  ,~~~~
\d{\bar\u}( t) = \lambda_{n}^\eps\l(\bar\xi(t), \bar\u(t)\r), 
\ee
with initial data 
$
\bar\xi(0)=\xi_s,~\bar\u(0)=\u_s
$
such that
\be
\tnv{v}{0}{T}^2+\tnw{w}{0}{T}^2\le \tilde C \eps^{2n}.\nn
\ee
The constant $\tilde C$ depends on  $F$ and $\xi_s$.
The parameters $\bar\xi,\bar\u$ describe
a fixed nontrivial perturbation of the uniform linear motion
as $\eps \to 0$ if condition \re{intro condition on nontrivial dynamic} is satisfied.  
\eth
\noindent

\boldmath\section{Construction of the Virtual Solitary Manifold }\unboldmath \la{Virtual Solitary Manifold}

\subsection{Iteration Scheme}\la{se:Implicit function theorem}
Let $\a,n\in\N$.
In this subsection, 
we establish
the iteration scheme presented in the introduction.
We implement the scheme for a general function $\ti F: (-1,1) \to H^{1,\a}(\R^2),~ \eps \mapsto \ti F(\eps)$ of class $C^{n}$ that satisfies $\ti F(0)=0$.
%
For being able to apply the implicit function theorem in the proof of existence of iterative solutions, we need to show that the corresponding linearizations of 
$(\hatt,\hatp,\lambda) \mapsto 
\ti {\cal G}_j^\eps(\hatt,\hatp,\lambda)
$ carried out at 
$(\hatt,\hatp,\lambda)=(0,0,0)$, $\eps=0$ 
%
are invertible ($\ti {\cal G}_j$ given by \re{successive eq G1}-\re{intro def tilde G}).
%
This is done in the following proposition, which is a main ingredient in the construction of the virtual solitary manifold. 
We start with a definition.
\bde\la{Def space Z}
\begin{itemize}
\item [(a)] $\ubar{ Z}^\a$ is the space $H^{2,\a}(\R^2) \oplus H^{1,\a}(\R^2)$ 
with the finite norm
$$
|z|_{\ubar{ Z}^\a} = |\v|_{H^{2,\a}(\R^2)}+ |\w|_{H^{1,\a}(\R^2)}.
$$

\item [(b)] $Z_n^\a=Z_n^\a(\u_*)$ is the space $\bigg\{ z =(\v,\w) \in C^n(I(\u_*), \ubar{ Z}^\a) : \Vert z \Vert_{Z_n^\a(\u_*)} <\infty \bigg\}\,$ 
with the finite norm
$$
\Vert z \Vert_{Z_n^\a(\u_*)} =\sup_{\u\in I(\u_*)} \l( \sum_{i=0}^n |\du^i y(u)|_{\ubar{ Z}^\a}\r).
$$
\item [(c)] We denote by
$t_{1}(\xi,\u,x):= \bma
\dxi\t_0(\xi,u,x)\\
\dxi\p_0(\xi,u,x)\\
\ema $
and by 
$t_{2}(\xi,\u,x):=\bma
\du\t_0(\xi,u,x)\\
\du\p_0(\xi,u,x)\\
\ema
$, where $\u\in(-1,1),~\xi,x\in\R $.
\end{itemize}
\ede

\noindent Recall that the spaces $\ubar{ Y}^\a$, $Y_n^\a(\u_*)$ were defined in \cref{ch: Main Results virtual}. We set $Y_n^\a:=Y_n^\a(\u_*)$.
\bpro \la{le invertibilityMxiCtwo alpha}
There exists $\underline{u}^\a>0$ such that the operator 
$
{\frak M}_n^\a:  Y_n^\a(u_*)  \to Z_n^\a(u_*),~\\
( \t,
\p,
\lambda) \mapsto
{\frak M}_n^\a
( \t,
\p,
\lambda ),
$
given by
\be
{\frak M}_n^\a
( \t,
\p,
  \lambda  )(\u)
=\bma
\u\dxi\t(u) -\p(u)  \\ 
-\dx^2\t(u) +\cos(\t_K(\Z))\t(u) +\u\dxi\p(u)  \\ 
\ema
+ \lambda(u)  
t_2(\xi,\u,x)
,\nn
\ee
is invertible if $u_*< \underline{u}^\a$.
\epro
\noindent
For proof see \cref{Preliminary Decompositions} and \cref{Invertibility of the Linearization}.
\noindent We formalize the iteration scheme in the following theorem. The maps $\tilde {\cal G}_j$ 
are defined on spaces of different regularity in $u$ such that the regularity of the spaces decreases with increasing $j$, which ensures well-definedness of $\tilde {\cal G}_j$.   
\bth\la{thimplicitfunctionIT1 alpha} 
Let $J=(-1,1)$, $u_*<\underline{u}^\a$
and let $\ti F: J \to H^{1,\a}(\R^2)\,, \eps \mapsto \ti F(\eps)$ be a $C^{n}$ function such that $\ti F(0)=0$. 
%
Let $\ti{\cal G}_1$ be given by
\be
{}&\ti{\cal G}_1: J  \times Y_{n+1}^\a (u_*) \to  Z_{n+1}^\a(u_*),~
(\eps,\hatt,\hatp,\lambda) \mapsto \tiGo(\hatt,\hatp,\lambda):=\Go(\t_0+\hatt,\p_0+\hatp,\lambda),\nn
\ee
where ${\cal G}_1$ is defined by \re{successive eq G1}. Then there exists $\eps^*>0$ and
a map
\be
{}&(-\eps^*,+\eps^*) \to Y_{n+1}^\a(u_*),~
\eps \mapsto (\hat\t_1^\eps,\hat\p_1^\eps,\lambda_{ 1}^\eps),\nn
\ee
of class $C^n$ such that
$
\ti{\cal G}_1^\eps(\hatt_1^\eps,\hatp_1^\eps,\lambda_{1}^\eps)=0.
$
Let $\ti{\cal G}_2$ be given by
\be
\ti{\cal G}_2: J \times Y_{n}^\a(u_*)  \to  Z_{n}^\a(u_*),~
 (\eps,\hatt,\hatp,\lambda) \mapsto \tiGt(\hatt,\hatp,\lambda):=\Gt(\t_0+\hatt,\p_0+\hatp,\lambda),\nn
\ee
where ${\cal G}_2$ is defined by \re{successive eq G2} with $(\t_1^\eps ,\p_1^\eps ,\lambda_{1}^\eps)
=(\t_0+\hatt_1^\eps ,\p_0+\hatp_1^\eps ,\lambda_{1}^\eps)$. Then there exists 
a map
\be
(-\eps^*,+\eps^*) \to Y_{n}^\a(u_*),~
\eps \mapsto (\hat\t_2^\eps,\hat\p_2^\eps,\lambda_{ 2}^\eps),\nn
\ee
of class $C^n$ such that
$
\ti{\cal G}_2^\eps(\hatt_2^\eps,\hatp_2^\eps,\lambda_{2}^\eps)=0\,.
$
This process can be continued successively to arrive at $\ti{\cal G}_n$ be given by
\be
\ti{\cal G}_n: J \times Y_{2}^\a(u_*)  \to  Z_{2}^\a(u_*),~
(\eps,\hatt,\hatp,\lambda) \mapsto \tiGn(\hatt,\hatp,\lambda):=\Gn(\t_0+\hatt,\p_0+\hatp,\lambda),\nn
\ee 
where ${\cal G}_n$ is defined by \re{successive eq Gn} with $(\t_{n-1}^\eps ,\p_{n-1}^\eps ,\lambda_{n-1}^\eps)
=(\t_0+\hatt_{n-1}^\eps ,\p_0+\hatp_{n-1}^\eps ,\lambda_{n-1}^\eps)$.
Ultimately there exists 
a map
\be
(-\eps^*,+\eps^*) \to Y_{2}^\a(u_*),~
\eps \mapsto (\hat\t_n^\eps,\hat\p_n^\eps,\lambda_{ n}^\eps),\nn
\ee 
of class $C^n$ such that
$
\ti{\cal G}_n^\eps(\hatt_n^\eps,\hatp_n^\eps,\lambda_{n}^\eps)=0\,
$
and we set $(\t_{n}^\eps ,\p_{n}^\eps ,\lambda_{n}^\eps)
=(\t_0+\hatt_{n}^\eps ,\p_0+\hatp_{n}^\eps ,\lambda_{n}^\eps)$.
\eth
\bpr
We skip $u_* $ in the notation.
Notice that
$
\ti{\cal G}_1^{0}(0,0,0)={\cal G}_1^{0}(\t_0,\p_0,0)=0\,.
$
The derivative of 
$
\ti{\cal G}_1: J \times Y_{n+1}^\a  \to  Z_{n+1}^\a
$ 
with respect to $(\hatt,\hatp,\lambda)$ evaluated at $(\eps,\hatt,\hatp,\lambda)=(0,0,0,0)$ 
is
$
{\frak M}_{n+1}^\a
$, which is invertible due to \cref{le invertibilityMxiCtwo alpha}. By the implicit function theorem there exists a 
$\eps_1^*>0$ and a map
$$
(-\eps_1^*,+\eps_1^*) \to Y_{n+1}^\a,
~\eps \mapsto (\hat\t_1^\eps,\hat\p_1^\eps,\lambda_{1}^\eps)
$$
of class $C^n$ such that
$
\tiGo(\hat\t_1^\eps,\hat\p_1^\eps,\lambda_{1}^\eps)=0.
$
We continue successively this process until we
obtain 
that the derivative of 
$
\ti{\cal G}_n: J \times Y_{2}^\a  \to  Z_{2}^\a
$ 
is
$
{\frak M}_{2}^\a.
$
This yields by using    
the same argument combined with $
\ti{\cal G}_n^{0}(0,0,0)
=0\,
$ that there exists
$\eps_n^*>0$ and a map
\be
(-\eps_n^*,+\eps_n^*) \to Y_{2}^\a,
~\eps \mapsto (\hat\t_n^\eps,\hat\p_n^\eps,\lambda_{ n}^\eps)\nn
\ee 
of class $C^n$ such that
$
\tiGn(\hat\t_n^\eps,\hat\p_n^\eps,\lambda_{ n}^\eps)=0\,.
$ 
We set $\eps^*=\min\{\eps_1^*,\eps_2^*,\ldots,\eps_n^*\}$.
\epr
\noindent
In the following theorem we state the properties of the $n$th iterative solution from \cref{thimplicitfunctionIT1 alpha}.  
\bth \la{thITrelations}
Let the assumptions of \cref{thimplicitfunctionIT1 alpha} hold. Then the following relations are satisfied.
\begin{itemize}
\item [(a)] $(\deps^j\t_{n-1}^{0},\deps^j\p_{n-1}^{0},\deps^j\lambda_{ {n-1}}^{0})=(\deps^j\t_n^{0},\deps^j\p_n^{0},\deps^j\lambda_{ n}^{0})$ for $j=0, \ldots, n-1$, where $n\ge 2$.
\item [(b)] $\forall \u\in I:~$      
$\ba\la{ITnquantitativ}
{}&
\u\dxi\l(\begin{matrix}
\t_n^\eps\\
\p_n^\eps\\
\end{matrix}\r)
-\l(\begin{matrix}
\p_n^\eps\\
[\t_n^\eps]\xx-\sin\t_n^\eps+\ti F(\eps)\\
\end{matrix}\r)
+\lambda_{ n}^\eps\du\bma
\t_n^\eps\\
\p_n^\eps\\
\ema +\Rn =0,\\
\ea
$
where\\
$
\ba
\Rn: =
\lambda_{ n}^\eps\du\bma
\sum_{i=0}^{n-1} \fr{\deps^i\t_{n}^0}{i!}\eps^i-\t_n^\eps\\
\sum_{i=0}^{n-1} \fr{\deps^i\p_{n}^0}{i!}\eps^i -\p_n^\eps\\
\ema.
\ea
$\\
The following rates of convergence hold. 
\be
{}&\l\Vert \bma
\sum_{i=0}^{n-1} \fr{\deps^i\t_{n}^0}{i!}\eps^i-\t_n^\eps\\
\sum_{i=0}^{n-1} \fr{\deps^i\p_{n}^0}{i!}\eps^i -\p_n^\eps\\
0\\
\ema \r\Vert_{Y_2^\a(\u_*)} = {\cal O}(\eps^n),~~~~~~
\l\Vert \bma
0\\
0\\
\lambda_{ n}^\eps\\
\ema \r\Vert_{Y_2^\a(\u_*)} = {\cal O}(\eps).\nn
\ee
\end{itemize}
\eth
\bre
The derivatives of the iterative solutions coincide at $0$ in the following way: 
$(\deps^j\t_{1}^{0},\deps^j\p_{1}^{0},\deps^j\lambda_{ {1}}^{0})
= (\deps^j\t_{2}^{0},\deps^j\p_{2}^{0},\deps^j\lambda_{{2}}^{0})$ for $j=0, 1;$
$(\deps^j\t_{2}^{0},\deps^j\p_{2}^{0},\deps^j\lambda_{ {2}}^{0})
= (\deps^j\t_{3}^{0},\deps^j\p_{3}^{0},\deps^j\lambda_{ {3}}^{0})$ for $j=0, 1,2 $ and so on up to the identities for the $n$th iterative solution stated in \cref{thITrelations} (a). 
\ere
\bpr[of \cref{thITrelations}]
(a) The claim can be proved by induction on $n$. We show the induction step.
Assume that the claim is true for all integers less than or equal to $n-2$.
Let $0\le j\le n-1$. The fact that 
the solutions
from \cref{thimplicitfunctionIT1 alpha} satisfy
$$
\forall u\in I(u_*):~~~(\t_j^\eps(u),\p_j^\eps(u),\lambda_{ j}^\eps(u))\in \ubar{ Y}^\a= H^{3,\a}(\R^2) \oplus H^{2,\a}(\R^2) \oplus H^{2,\a}(\R)\,
$$
and that the injections
$H^1(\R) \subset L^\infty(\R)\,,
H^2(\R^2) \subset L^\infty(\R^2)\,$
are continuous \cite[Corollary 9.13]{Brezis} yields the justification for using the Leibniz's formula and Fa\`a di Bruno's formula.
Thus we obtain for the
$j$-th derivatives with respect to $\eps$, evaluated at
$\eps=0$:
\be\la{ITnz}
0={}&\deps^j{\cal G}_n^0(\t_n^0,\p_n^0,\lambda_{ n}^0)\\
={}&
\bma
\u\dxi\deps^j\t_n^0-\deps^j\p_n^0\\
\u\dxi\deps^j\p_n^0-\dx^2\deps^j\t_n^0
\ema
+\bma
0\\
 \sum_{I_j} \fr {j!}{l_1!l_2!\ldots l_j!} \deps^l \sin (\t_n^0) \l( \fr {\deps^1 \t_n^0}{1!} \r)^{l_1}
\l( \fr {\deps^2 \t_n^0}{2!}  \r)^{l_2}\ldots
\l( \fr {\deps^j \t_n^0}{j!}  \r)^{l_j}
\ema
\nn
\\
{}&
-\bma
0\\
 \deps^j  \ti F(0)
\ema
+ \bma
\sum_{i=0}^j \binom{j}{i}   \deps^{j-i}\lambda_{ n}^0  {\du\deps^i\t_{n-1}^0} \\
\sum_{i=0}^j \binom{j}{i} \deps^{j-i}\lambda_{ n}^0  {\du\deps^i\p_{n-1}^0} \\
\ema,\nn
\ee

\be\la{ITnmoz}
0={}&\deps^j{\cal G}_{n-1}^0(\t_{n-1}^0,\p_{n-1}^0,\lambda_{ n-1}^0)\\
={}&
\bma
\u\dxi\deps^j\t_{n-1}^0-\deps^j\p_{n-1}^0\\
\u\dxi\deps^j\p_{n-1}^0-\dx^2\deps^j\t_{n-1}^0
\ema
+\bma
0\\
 \sum_{I_j} \fr {j!}{l_1!l_2!\ldots l_j!} \deps^l \sin (\t_{n-1}^0) \l( \fr {\deps^1 \t_{n-1}^0}{1!} \r)^{l_1}
\l( \fr {\deps^2 \t_{n-1}^0}{2!}  \r)^{l_2}\ldots
\l( \fr {\deps^j \t_{n-1}^0}{j!}  \r)^{l_j}
\ema
\nn\\
{}&
-\bma
0\\
 \deps^j \ti F(0)
\ema
+ \bma
\sum_{i=0}^j \binom{j}{i}  \deps^{j-i} \lambda_{ n-1}^0 {\du\deps^i\t_{n-2}^0}\\
\sum_{i=0}^j \binom{j}{i}  \deps^{j-i} \lambda_{ n-1}^0 {\du\deps^i\p_{n-2}^0} \\
\ema,\nn
\ee
where $l= l_1+l_2+\ldots +l_j$ and the sum is taken over all $l_1,l_2,\ldots ,l_j$ for which $l_1+ 2l_2+\ldots +jl_j=j$.
Subtracting \re{ITnmoz} from \re{ITnz} yields the claim due to \cref{le invertibilityMxiCtwo alpha}.\\
(b) follows from (a), \cref{thimplicitfunctionIT1 alpha} and Taylor's formula.
\epr
\subsection{Virtual Solitary Manifold}
From now on we set $\a:=1$. 
In this subsection, we apply \cref{thimplicitfunctionIT1 alpha} on a specific $\ti F$ and define the virtual solitary manifold by the solution obtained in the $n$th iteration.  
\bde\la{de cutoff}
Let $F,\xi_s$ be from \cref{maintheorem specialF kderivativesvanish virtual manifold} and $\Xi:=\Xi(\xi_s):=  |\xi_s|+3$.
We set
$\ti F(\eps,\xi,x):=\F \chi(\xi),$ 
where $\chi$ is a smooth cutoff function
with $\chi (\xi)=1$ for $|\xi|\le \Xi$ and $\chi (\xi)=0$ for $|\xi|\ge \Xi+1$.
\ede
\noindent
The next lemma follows immediately from \cref{thITrelations} and from the assumptions on $F$ in \cref{maintheorem specialF kderivativesvanish virtual manifold}.
\ble \la{le: Ftilde vs F}
Let $F$ be from \cref{maintheorem specialF kderivativesvanish virtual manifold}, $ \ti F$ from \cref{de cutoff}.
Then it holds that
\begin{itemize}
\item[(a)] $\forall ~(\eps,\xi,x) \in (-1,1)\times \l[-\Xi,\Xi \r] \times \R: \ti F(\eps,\xi,x)=\F$.
\item[(b)] 
$\ti F \in C^{n}((-1,1),H^{1,1}(\R^2))$ and  
$\deps^l \ti F(0,\cdot,\cdot)=0$ for $0\le l\le k 
.$
\item[(c)]
$
\ba
\l\Vert \bma
0\\
0\\
\lambda_{n}^\eps\\
\ema \r\Vert_{Y_2^1(\u_*)} = {\cal O}(\eps^{k+1}).
\ea
$
\end{itemize}
\ele
\ble\la{le transestimate}
Let $v\in H^1(\R^2)$. Then there exists $ b>0$ such that 
\be
\forall ~\xi\in \R ~\nltwox{v(\xi,x)}\le b\nhonertwoxix{v(\xi,x)}.\nn
\ee
\ele

\bpr
This follows from applying Morrey's embedding Theorem to the variable $\xi$.
\epr
\noindent
We solve iteratively the equations in 
\cref{thimplicitfunctionIT1 alpha} with the specific $\ti F(\eps,\xi,x):=\F \chi(k,\xi)$ from \cref{de cutoff}
and define by the $n$th solution $(\t_n^\eps,\p_n^\eps,\lambda_n^\eps)$ the virtual solitary manifold ${\cal S}_n^\eps$.
We utilize the truncated version of $F$ rather than $F$ itself in order to make sure that the maps $\ti {\cal G}_j$ in \cref{thimplicitfunctionIT1 alpha} are well defined. 
\cref{thimplicitfunctionIT1 alpha} is applicable to $\ti F$ 
due to \cref{le: Ftilde vs F}. 
\bde\la{def:virtual solitary manifold}
Let $\underline{u}^1$ be from \cref{le invertibilityMxiCtwo alpha} (case $\a=1$). We fix a specific $u_*$ such that $0<u_*<\underline{u}^1$. 
Let $\ti F$ be from \cref{de cutoff}. 
Let $\eps^*>0$ be the constant and let $(\t_n^\eps,\p_n^\eps,\lambda_n^\eps)$ be the nth solution obtained from application of \cref{thimplicitfunctionIT1 alpha} to $\ti F$.
We set
\be
{\cal S}_n^\eps:=\l\{ \bma
\t_n^\eps(\xi,\u,\cdot)\\
\p_n^\eps(\xi,\u,\cdot)
\ema~:~\u\in(-u_*,u_*),~\xi\in\R
\r\}~\text{for}~\eps \in (-\eps^*,\eps^* ),\nn
\ee
and call ${\cal S}_n^\eps$ the virtual solitary manifold.
\ede
\bre 
\mbox{}
\begin{itemize}
\item[(a)] 
From now on we 
denote by $(\t_{n}^\eps ,\p_{n}^\eps ,\lambda_{n}^\eps)
=(\t_0+\hatt_{n}^\eps ,\p_0+\hatp_{n}^\eps ,\lambda_{n}^\eps)$ always the $n$th solution utilized in \cref{def:virtual solitary manifold}.
\item[(b)] The vectors
\la{de: virtual tangent vectors}
$$
t_{1,n}^\eps(\xi,\u,x):= \bma
\dxi\ttn\\
\dxi\ptn\\
\ema ~~~~\text{and}~~~~ t_{2,n}^\eps(\xi,\u,x):=\bma
\du\ttn\\
\du\ptn\\
\ema
$$
are tangent vectors of the manifold ${\cal S}_n^\eps$ at the point $( \t_n^\eps(\xi,\u,\cdot), \p_n^\eps(\xi,\u,\cdot) )$ and form a basis of the tangent space at this point.  
\end{itemize} 
\ere
\section{Symplectic Orthogonal Decomposition}\la{ch:Symplectic Orthogonal Decomposition}
Let from now on $u_*,\eps^*$ be always from \cref{def:virtual solitary manifold} and let $U$ be fixed such that $0<U<u_*$. 
We consider $V(l,U,u_*),~\Sigma(l,U,u_*)$ introduced in \cref{def:PartfourMainResult} (e)
and the function $
{\cal N}^\eps : L^\infty(\R) \times L^2(\R)\times \Sigma(2,U,u_*)\to\R^2 
$ defined as in \re{orthogonalitycondIntro} for these specific $u_*,U$. For simplicity of further notation we set $\Sigma(l,U,u_*)=\Sigma(l)$ and $V(l,U,u_*)=V(l)$.
\\
In this chapter we will choose $\eps_0$  sufficiently small and consider $\eps\in(0,\eps_0]$. 
We show that if $(\t,\p)\in  L^\infty(\R)\oplus L^2(\R)$ is close enough (in the $L^\infty(\R)\oplus L^2(\R)$ norm) to the region 
\be
{\cal S}_n^\eps(U):=\l\{ \bma
\t_n^\eps(\xi,\u,\cdot)\\
\p_n^\eps(\xi,\u,\cdot)
\ema~:~(\xi,\u)\in \Sigma(4) 
\r\},\nn
\ee
of the virtual solitary manifold 
${\cal S}_n^\eps$, then there exists a unique $(\xi,\u)\in \Sigma(2)$ such that we are able to decompose the solution
$$
\bma
\t(\cdot)\\
\p(\cdot)
\ema
= \bma
\t_n^\eps(\xi,\u,\cdot)\\
\p_n^\eps(\xi,\u,\cdot)
\ema+
\bma
\v(\cdot)\\
\w(\cdot)
\ema,
$$
in a point on the virtual solitary manifold $(\t_n^\eps(\xi ,\u ,\cdot),\p_n^\eps(\xi ,\u ,\cdot))$  and a transversal component $(v( \cdot),w( \cdot))$, which
is symplectic orthogonal to the tangent vectors 
$t_{1,n}^\eps(\xi,\u,\cdot)$
and
$t_{2,n}^\eps(\xi,\u,\cdot)$
at the corresponding point of the manifold ${\cal S}_n^\eps$,
i.e., the orthogonality condition 
$$
{\cal N}^\eps(\t,\p,\xi,\u)=0\,
$$
is satisfied.
We prove that the symplectic decomposition is  possible in a small uniform distance to the manifold ${\cal S}_n^\eps$, where the distance might depend on $\eps_0$ but does not depend on $\eps$. 
\bre
In \cref{thimplicitfunctionIT1 alpha} we have solved the equations defining $(\t_n^\eps,\p_n^\eps,\lambda_{n}^\eps)$ in weighted spaces. One of the reasons for working in weighted Sobolev spaces was to make sure that ${\cal N}^\eps : L^\infty(\R) \times L^2(\R)\times \Sigma(2)\to\R^2$ is 
well defined.  
\ere
\noindent
We start with a definition and some elementary lemmas which will be used later.
\bde\la{decomposition of virtual manifold}
Let $\eps \in(0,\eps^*)$. We set
\begin{itemize}
\item[(a)] $  m:= \int  [\t_K'(Z)]^2\,dZ $,
\item[(b)] $\ba m_n^\eps(\xi,\u):=\int -\dxi\ptn \du\ttn+\dxi\ttn \du\ptn \,dx,\ea$
\item[(c)]
$\ba
k_n^\eps(\xi,\u)  :={}& \int   - \dxi\p_0(\xi,\u,x)\du\hat\t_n^\eps(\xi,\u,x)
- \du \t_0(\xi,\u,x) \dxi\hat\p_n^\eps(\xi,\u,x)  \\
{}& +\dxi \t_0(\xi,\u,x) \du\hat\p_n^\eps(\xi,\u,x)  
  +\du \p_0(\xi,\u,x) \dxi\hat\t_n^\eps(\xi,\u,x) \\
{}&-\dxi\hat\p_n^\eps(\xi,\u,x)\du\hat\t_n^\eps(\xi,\u,x)+\dxi\hat\t_n^\eps(\xi,\u,x)\du\hat\p_n^\eps(\xi,\u,x)  
\,dx.
\ea$
\end{itemize}
\ede
\noindent A straight forward computation yields the following lemma.
\ble \la{le:mepsidentity}
Let $\eps \in(0,\eps^*)$. It holds that
\be
\forall ~(\xi,u)\in\R\times [-U- V(2),U+ V(2)]:~~~m_n^\eps(\xi,\u)
 = {}& \g^3(\u)m +k_n^\eps(\xi,\u).\nn
\ee 
\ele
\ble\la{le controlvirtualtvproduct}
Let $\eps_0>0$ be sufficiently small. 
There exist constants
$ c=c(U)>0,~C=C(U)>0,$ 
such that $\forall \eps\in(0,\eps_0],~(\xi,u)\in\R\times [-U- V(2),U+ V(2)]:~$ 
\be
 c \le \fr {\g^3(\u)m} 2\le m^\eps_n(\xi,\u) \le  2{\g^3(\u)m} \le C.\nn
\ee
\ele
\bpr
Using \cref{le transestimate} and 
continuity of 
$
\eps \mapsto (\hat\t_n^\eps,\hat\p_n^\eps,\lambda_{n}^\eps)$
(see \cref{thimplicitfunctionIT1 alpha}) we obtain for sufficiently small $\eps_0$:  
$
\forall \eps\in(0,\eps_0],~(\xi,u)\in\R\times [-U- V(2),U+ V(2)]:|k^\eps_n(\xi,\u)|<m/ 2,
$
which implies the claim.
\epr
\noindent 
The next lemma provides that the symplectic decomposition described above is possible. In the proof we will take derivatives of $(\t_n^\eps,\p_n^\eps)$ up to second order with respect to $\xi$ and $u$. This was the reason for solving, in \cref{Virtual Solitary Manifold}, the equations defining $(\t_n^\eps,\p_n^\eps,\lambda_{n}^\eps)$ in spaces of higher regularity in $\xi$ and $\u$.  
\ble\la{le uniform decomposition virtual} 
Let $\eps_0>0$ be sufficiently small.
Let
\be
{\cal O}={\cal O}^\eps_{U,p}=\Big\{(\t,\p)\in L^\infty(\R)\times L^2(\R): \inf_{{(\xi,\u)}\in  \Sigma(4)} \l| \bma \t(\cdot)\\ \p (\cdot)\ema - \bma\t_n^\eps(\xi,\u,\cdot)\\ \p_n^\eps(\xi,\u,\cdot)\ema\r|_{L^\infty(\R)\oplus L^2(\R)} <p\Big\}.\nn
\ee
There exists $r>0$ such that if $\eps\in(0,\eps_0]$
and $p \le r$ then for any $(\t, \p) \in {\cal O}^\eps_{U,p}$ there exists a unique 
$(\xi,\u)\in \Sigma(2)$ such that 
\be
{\cal N}^\eps(\t,\p,\xi,\u)=0\,\nn
\ee
and the map 
$(\t,\p) \mapsto (\xi(\t,\p) ,\u(\t,\p) )$ 
is in $C^1({\cal O}^\eps_{U,p}, \Sigma(2))$.
\ele
\bpr
Let $\eps_0\in(0,\eps^*)$.
We will specify $\eps_0$ later in this proof.
Let $\eps\in(0,\eps_0]$. Notice that the map 
$
\eps \mapsto (\hat\t_n^\eps,\hat\p_n^\eps,\lambda_{n}^\eps)$
from \cref{thimplicitfunctionIT1 alpha}
is continuous and it holds 
$\Sigma(4)
\subset \Sigma(3)\subset \Sigma(2)$.
Consider $(\xi_0, u_0)\in\Sigma(3)$.
\cref{le:mepsidentity} yields that
\be\la{DxiuNeps}
{}&D_{\xi,\u} {\cal N}^\eps(\t_n^\eps(\xi_0,\u_0,\cdot),\p_n^\eps(\xi_0,\u_0,\cdot),\xi_0,\u_0)
=
\l(\g^3(\u_0)m+ k_n^\eps(\xi_0,\u_0) \r)\bma
0
& 1 \\
-1
& 0
\ema.
\ee
Using \cref{le transestimate} we obtain for sufficiently small $\eps_0$ for all 
$\eps\in(0,\eps_0]$:
$
\la{deltaeps}
|k_n^\eps(\xi_0,\u_0)|\le \fr m 2
$   
and thus
\be\la{detnotzero}
\det D_{\xi,\u} {\cal N}^\eps(\t^\eps(\xi_0,\u_0,\cdot),\p^\eps(\xi_0,\u_0,\cdot),\xi_0,\u_0)\not=0.
\ee
We prove that there exist $r>0,{\bar\delta}>0,\eps_0>0$
such that 
$\forall\eps\in(0,\eps_0], (\xi_0,\u_0)\in \Sigma(3)$ there exist balls 
$
B_r(\t_n^\eps(\xi_0,\u_0,\cdot),\p_n^\eps(\xi_0,\u_0,\cdot)) \subset L^\infty(\R)\oplus L^2(\R)\,$,
$B_{\bar\delta}(\xi_0,\u_0)\subset \Sigma(2)\,,$
and a map
$$
T_{\xi_0,\u_0}^\eps:  B_r(\t_n^\eps(\xi_0,\u_0,\cdot),\p_n^\eps(\xi_0,\u_0,\cdot)) \to B_{\bar\delta}(\xi_0,\u_0)
$$
such that 
$
{\cal N}^\eps(\t,\p,T_{\xi_0,\u_0}^\eps(\t,\p))=0
$ 
on
$B_r(\t_n^\eps(\xi_0,\u_0,\cdot),\p_n^\eps(\xi_0,\u_0,\cdot))$. Therefore we refer to \cite[Theorem 15.1]{Deimling} and check their proof of the implicit function theorem, whereas we show that $r$ and ${\bar\delta}$ do not depend on $\eps$ and on $(\xi_0,\u_0)$. 
We introduce  
$$
\bar{\cal  N}_{\xi_0,\u_0}^\eps(\t,\p,\xi,\u):={\cal  N}^\eps(\t(\cdot)+\t_n^\eps(\xi_0,\u_0,\cdot),\p(\cdot)+\p_n^\eps(\xi_0,\u_0,\cdot),\xi+\xi_0,\u+\u_0).$$
Notice that
$
\bar{\cal N}_{\xi_0,\u_0}^\eps(0,0,0,0)=(0,0).
$
We set
$
K_{\xi_0,\u_0}^\eps :=
D_{(\xi,\u)} \bar{\cal N}_{\xi_0,\u_0}^\eps(0,0,0,0)
$
and
\be
 S_{\xi_0,\u_0}^\eps(\t,\p,\xi,\u):={}& \l[K_{\xi_0,\u_0}^\eps\r]^{-1}\bar{\cal N}_{\xi_0,\u_0}^\eps(\t,\p,\xi,\u)- I (\xi,\u),\nn
\ee
which is well defined due to \re{detnotzero}.
Due to \cref{le transestimate} it holds for a  sufficiently small $\eps_0$ that
\be
\begin{split}
{}&\forall \eps\in(0,\eps_0], (\xi,u)\in\R\times [-U- V(2),U+ V(2)],~\beta_1+\beta_2\le 2,~ p=1,2:\\
{}&\l|\dxi^{\beta_1}\du^{\beta_2}\t^\eps(\xi,\u,x)\r|_{L^p_x(\R)}\le B,~~~~~~
\l|\dxi^{\beta_1}\du^{\beta_2}\p^\eps(\xi,\u,x)\r|_{L^p_x(\R)}\le B.
\end{split}\nn
\ee
\noindent
In this proof we denote by $\Vert \cdot\Vert$ the maximum row sum norm of a $2\times 2$ matrix induced by the maximum norm $|\cdot|_{\infty}$ in $\R^2$.
We claim that $\exists k\in (0,1),{\bar\delta} >0,\eps_0>0~~\forall \eps\in(0,\eps_0] , (\xi_0,\u_0)\in \Sigma(3)  
$
$\forall \l((\t,\p),(\xi,\u)\r)\in B_{\bar\delta}(0)\times B_{\bar\delta}(0):~\Vert D_{(\xi,\u)}S_{\xi_0,\u_0}^\eps(\t,\p,\xi,\u) \Vert \le k <1$.
Due to \re{DxiuNeps} it holds that
\be
{}&D_{(\xi,\u)}S_{\xi_0,\u_0}^\eps(\t,\p,\xi,\u)\nn\\
={}& \fr 1{\gz^3 m+ k^\eps_n(\xi_0,\u_0)} 
\bma
-\dxi\bar{\cal  N}_{\xi_0,\u_0,2}^{\eps}(\t,\p,\xi,\u)  & -\du\bar{\cal  N}_{\xi_0,\u_0,2}^{\eps}(\t,\p,\xi,\u)\\
 \dxi\bar{\cal  N}_{\xi_0,\u_0,1}^{\eps}(\t,\p,\xi,\u) & \du\bar{\cal  N}_{\xi_0,\u_0,1}^{\eps}(\t,\p,\xi,\u)
\ema
- \bma
1 & 0\\
0 & 1
\ema.\nn
\ee
The claim follows by using \cref{le transestimate}, \cref{le controlvirtualtvproduct} and estimating each entry of $D_{(\xi,\u)}S_{\xi_0,\u_0}^\eps(\t,\p,\xi,\u)$, for instance:
\be
{}&|-\fr 1 {m^\eps_n(\xi_0,\u_0)}\dxi\bar{\cal  N}_{\xi_0,\u_0,2}^{\eps}(\t,\p,\xi,\u)-1|\nn\\
\le{}&
\fr 1 {|m^\eps_n(\xi_0,\u_0)|}\Big(
|\dxi\du\ptbn|_{L^2_x(\R)}|\hat\t_n^\eps(\xi_0,\u_0,x)-\hat\t_n^\eps(\bar\xi,\bar\u,x)|_{L^2_x(\R)}\nn\\
{}&+|\dxi\du\ptbn|_{L^1_x(\R)}|\t(x)+\t_0(\xi_0,\u_0,x)-\t_0(\bar\xi,\bar\u,x)|_{L^\infty_x(\R)}\nn\\
{}&+|\dxi\du\ttbn|_{L^2_x(\R)}|\p(x)+\hat\p_n^\eps(\xi_0,\u_0,x)-\hat\p_n^\eps(\bar\xi,\bar\u,x)|_{L^2_x(\R)}\nn\\
{}&+|\dxi\du\ttbn|_{L^1_x(\R)}|\p_0(\xi_0,\u_0,x)-\p_0(\bar\xi,\bar\u,x)|_{L^\infty_x(\R)}
+|m^\eps_n(\bar\xi,\bar\u)-m^\eps_n(\xi_0,\u_0)|\Big).\nn
\ee
Similarly as above one shows that $\exists r\le {\bar\delta},\eps_0>0~~ \forall \eps\in(0,\eps_0] , (\xi_0,\u_0)\in \Sigma(3) ~  
$
$\forall (\t,\p)\in B_r(0):~|S_{\xi_0,\u_0}^\eps(\t,\p,0,0) |_{\infty}< {\bar\delta}(1-k)$,
which completes the proof.
\epr
\section{Existence of Dynamics and the Orthogonal Component}\la{local solution F virtual}
We argue similar to \cite[Proof of theorem 2.1]{Stuart1}. Let $\eps_0$ be from \cref{le uniform decomposition virtual} and $\eps\in(0,\eps_0]$. In order to make use of existence theory we  
consider the problem
\be 
\bma
\bar\v(0,x)\\
\bar\w(0,x)
\ema {}&=
\bma
\t(0,x)-\t_n^\eps(\xi_s,\u_s,x)\\
\p(0,x)-\p_n^\eps(\xi_s,\u_s,x)
\ema\la{vweqn ID virtual F},\\
\dt\bma
\bar\v(t,x)\\
\bar\w(t,x)
\ema
{}&=\bma
\bar\w(t,x)-\p_n^\eps(\xi_s,\u_s,x)\\
[\bar\v(t,x)+\t_n^\eps(\xi_s,\u_s,x)]_{xx}-\sin (\bar\v(t,x)+\t_n^\eps(\xi_s,\u_s,x))+ \F \la{vweqn virtual F} \\
\ema.
\ee
By \cite[Theorem VIII 2.1, Theorem VIII 3.2
]{Martin} there exists a local solution (see also \cite[Proof of theorem 2.1]{Stuart1}, \cite[p.434
]{Stuart2}), where
$$
(\bar \v,\bar\w)\in C^1([0,T_{loc}], H^1(\R) \oplus L^2(\R)).
$$ 
$(\t,\p)$ given by $\t(t,x)=\bar\v(t,x)+\t_n^\eps(\xi_s,\u_s,x)$ and $\p(t,x)=\bar\w(t,x)+\p_n^\eps(\xi_s,\u_s,x) $ solves obviously locally the Cauchy problem \re{SGE1} and $(\t,\p) \in C^1([0,T_{loc}], L^\infty(\R) \oplus L^2(\R))$ due to Morrey's embedding theorem. 
We are going to obtain a bound in \cref{aprioriestimate F virtual manifold} which will imply that the local solutions are indeed continuable. 
So from now we assume that $(\bar\v,\bar\w) \in C^1([0,\overline T], H^1(\R) \oplus L^2(\R))$ is a solution of \re{vweqn ID virtual F}-\re{vweqn virtual F}	and $(\t,\p)$ is a solution of \re{SGE1} such that $(\t,\p) \in C^1([0,\overline T], L^\infty(\R) \oplus L^2(\R))$, where $\overline T>0$.
\\
In the following we define, similar to $\Sigma(l,U,u_*)$, a new parameter area, where the parameter $\xi$ is bounded. 
\bde \la{de strip orthogonality condition}
We set for $u_*, U$ as in \cref{ch:Symplectic Orthogonal Decomposition} and $\Xi$ from \cref{de cutoff} 
$$\Sigma(l,\Xi):= \Sigma(l,U,u_*,\Xi):=\Big\{  
(\xi,\u)\in 
 (-\Xi+1 - V(l) ,\Xi-1 + V(l)) \times (-U- V(l),U+ V(l)) \Big\}.$$
\ede
\noindent Given $(\t,\p)$ we choose the parameters $(\xi(t),\u(t))$ according to \cref{le uniform decomposition virtual}  and define $(v,w)$ as follows:
\be
v(t,x){}&=\t(t,x)- \tttn\la{decomposition1 F virtual},\\
w(t,x){}&=\p(t,x)- \pttn\la{decomposition2 F virtual}.
\ee
$(v(t,x),w(t,x))$ is well defined for $t \ge 0$ so small that 
$
\nlinf{v(t)}+\nltwo{w(t)}\le r\,
$
and
$
(\xi(t),\u(t))\in \Sigma(4,\Xi),
$
where $r$ is from \cref{le uniform decomposition virtual}. 
We formalize this in the following definition.
%
\bde\la{de: exittime decomposition virtual case}
Let $r$ be from \cref{le uniform decomposition virtual}. $t^*$ is the "exit time"
\be
t^*:= \sup\Big\{{}& T>0:\tnvinf{v}{0}{t}+\tnw{w}{0}{t}\le {r },
 (\xi(t),\u(t))\in \Sigma(4,\Xi),~0\le t \le T\Big\}.\nn
\ee
\ede
\noindent
Notice that $(\xi_s,\u_s)=(\xi(0),\u(0))\in \Sigma(4,\Xi) $. 
The transversal component $(v(t,x),w(t,x))$ is well defined for $0\le t \le t^*$
since, among others, we choose $\eps$ such that  
$\eps\in(0,\eps_0]$ with $\eps_0$ from \cref{le uniform decomposition virtual} and the initial data such that $\nlinf{v(0)}+\nltwo{w(0)}\le \fr r 2$ with $(\v(0),\w(0))$ given by \re{SGE1}.
\ble
Let $T=\min\{ t^*, \overline T\}$ and let $(v,w)$ be defined by \re{decomposition1 F virtual}-\re{decomposition2 F virtual}. Then
$(\v, \w) \in C^1([0,T], H^1(\R) \oplus L^2(\R))$. 

\ele
\bpr
This follows by using \re{decomposition1 F virtual}-\re{decomposition2 F virtual} 
and the fact that $
(\bar \v,\bar\w)\in C^1([0,T ], H^1(\R) \oplus L^2(\R)) 
$, since the difference $( \t_K(\g(u_0)(\cdot-\xi_0))-\t_K(\g(\bar \u)(\cdot-\bar\xi)) )$ is in $ L^2(\R)$ for all $(\xi_0,u_0),(\bar\xi,\bar u)\in \R\times (-1,1)$.
\epr
\noindent
In the following lemma we point out the relation between $F$ and $(\t_n^\eps,\p_n^\eps,\lambda_n^\eps)$. Notice that there appears $F$ instead of $\ti F$ in the equation above. 
\noindent
Moreover, we state the rates of convergence of $\Rn (\xi,u,\cdot)$ and $\lambda_{ n}^\eps(\xi,u)$ which will be needed in the proof of the modulation equations for the parameters $(\xi(t),\u(t)) $ in the next section and in the proof of the main result in \cref{aprioriestimate F virtual manifold}. 
\ble\la{le: remainder uniformly bound}
It holds that for a.e. $(\xi,\u,x)\in \Sigma(4,\Xi) \times \R$
\be
{}&
\u\dxi\l(\begin{matrix}
\t_n^\eps(\xi,u,x)\\
\p_n^\eps(\xi,u,x)\\
\end{matrix}\r)
-\l(\begin{matrix}
\p_n^\eps(\xi,u,x)\\
\dx^2\t_n^\eps(\xi,u,x)-\sin\t_n^\eps(\xi,u,x)+F(\eps,x)\\
\end{matrix}\r)\nn\\
{}&+\lambda_{ n}^\eps(\xi,u)\du\bma
\t_n^\eps(\xi,u,x)\\
\p_n^\eps(\xi,u,x)\\
\ema
 + \Rn (\xi,u,x)
=0\nn
\ee
and 
$
\nltwo{[\Rn(\xi,\u,\cdot)]_1} =  {\cal O}(\eps^{n+k+1}),
$   
$
\nltwo{[\Rn(\xi,\u,\cdot)]_2}  =  {\cal O}(\eps^{n+k+1}),
$
$
|{\lambda_{ n}^\eps(\xi,u)}| =  {\cal O}(\eps^{k+1}), 
$
$
|{\partial_1\lambda_{ n}^\eps(\xi,u)}| =  {\cal O}(\eps^{k+1}),
$
$|{\partial_2\lambda_{ n}^\eps(\xi,u)}|  =  {\cal O}(\eps^{k+1})$
uniformly in $(\xi,\u)\in \Sigma(4,\Xi)$.
\ele
\bpr
The first identity follows due to \cref{ITnquantitativ} and \cref{le: Ftilde vs F}.
Using \cref{ITnquantitativ}, \cref{le transestimate} and  Morrey's embedding theorem we obtain for all $(\xi,\u)\in \Sigma(4,\Xi)$:  
\be
\nltwo{[\Rn(\xi,\u,\cdot)]_1}
\le 
c \l\Vert \bma
0\\
0\\
\lambda_{ n}^\eps\\
\ema \r\Vert_{Y_2^\a(\u_*)}
\l\Vert \bma
\sum_{i=0}^{n-1} \fr{\deps^i\t_{n}^0}{i!}\cdot\eps^i-\t_n^\eps\\
\sum_{i=0}^{n-1} \fr{\deps^i\p_{n}^0}{i!}\cdot\eps^i -\p_n^\eps\\
0\\
\ema \r\Vert_{Y_2^\a(\u_*)}
=
{\cal O}(\eps^{n+k+1})\,\nn
\ee
and
\be
 |{\partial_1\lambda_{ n}^\eps(\xi,u)}|
\le 
%
c \l\Vert \bma
0\\
0\\
\lambda_{ n}^\eps\\
\ema \r\Vert_{Y_2^\a(\u_*)}
%
=
{\cal O}(\eps^{k+1}).\nn
\ee
The other cases can be treated analogously.
\epr
\noindent
We compute the time derivatives of $v$ and $w$, which will be needed in the following sections.
\ble\la{dvdw epstwo virtual}
Let $(v,w)$ be given by \re{decomposition1 F virtual}-\re{decomposition2 F virtual}. Then it holds
\be
\d \v(x) =  {}&\w(x) - \d\xi\dxi\ttn -\d\u\du\ttn \nn\\
{}&+\u\dxi\ttn  
+\luen(\xi,\u)\du\ttn
+\Rtno ,\nn\\ 
\d \w(x)   ={}&\dx^2\v(x) -\cos\ttn \v(x) +\fr{\sin\ttn\v^2(x) }{2} +\tilde R(\v)(x) 
+\u\dxi\ptn   \nn\\
{}&
+\luen(\xi,\u)\du\ptn +\Rtnt
- \d\xi\dxi\ptn -\d\u\du\ptn , \nn
\ee
for times $t \in [0,t^*]$, where $\tilde R(\v)={\cal O}(\nhonex{v}^3)$ and $\Rn (\xi,\u,x)$ is from \cref{thITrelations} (b).
\ele 
\bpr
By taking the time derivatives of $(v,w)$ and using \cref{le: remainder uniformly bound}, \re{SGE1} we obtain
\be
\d \v(x)  ={}& \w(x) +\ptn 
- \d\xi\dxi\ttn -\d\u\du\ttn \nn\\
= {}&\w(x) - \d\xi\dxi\ttn -\d\u\du\ttn 
\nn\\
{}&
+\u\dxi\ttn  
+\luen(\xi,\u)\du\ttn
+\Rtno\,\nn
\ee
and
\be
\d \w(x)  ={}& \dx^2\t(x) -\sin\t(x)+\F
- \d\xi\dxi\ptn -\d\u\du\ptn \nn\\
={}&
\dx^2\ttn+\dx^2\v(x) -\sin\ttn
-\cos\ttn\v(x)\nn\\
{}&+\fr{\sin\ttn\v^2(x)}{2}+\tilde R(\v)(x)+\F
- \d\xi\dxi\ptn -\d\u\du\ptn\nn\\
={}&\dx^2\v(x)-\cos\ttn \v(x)+\fr{\sin\ttn\v^2(x)}{2}+\tilde R(\v)(x)
+\u\dxi\ptn\nn\\
{}&   
+\luen(\xi,\u)\du\ptn +\Rtnt
- \d\xi\dxi\ptn -\d\u\du\ptn,\nn
\ee
where we have expanded the term $\sin(\ttn+\v(x))$.
\epr
\section{Modulation Equations}\la{ch Modulation Equations}
In the following lemma we derive modulation equations for the parameters $(\xi(t),\u(t))$.
\ble\la{le Modulation Equations}
There exists an $\eps_0>0$ such that the following statement holds. Let $\eps\in(0,\eps_0]$ and let $(v,w)$ be given by \re{decomposition1 F virtual}-\re{decomposition2 F virtual} 
with $(\xi,\u)$ obtained from \cref{le uniform decomposition virtual}. Let
$$ \tnv{v}{0}{t^*},\tnw{w}{0}{t^*}\le \eps_0. $$
Then it holds for $t \in [0,t^*]$ that
\be
|\d\xi(t)-u(t)|  \le{}&  C[\nhone{v(t)}+\nltwo{w(t)}] \eps^{k+1} + C\nhone{v(t)}^2 +C \eps^{n+k+1} ,\nn\\
	|\d\u(t)-\luen(\xi(t),\u(t))| \le{}&  C[\nhone{v(t)}+\nltwo{w(t)}] \eps^{k+1} + C\nhone{v(t)}^2 +C \eps^{n+k+1} \nn,
\ee
where C depends on $F$ and $\xi_s$.
\ele
\bpr
The technique we use is similar to that in the proof of  \cite[Lemma 6.2]{ImaykinKomechVainberg}.
Let $\eps_0\in(0,\eps^*)$ with $\eps^*$ from \cref{def:virtual solitary manifold} and let $\eps\in(0,\eps_0)$. Further in the proof we will make some more assumptions on $\eps_0$. 
We start with some definitions and set
\be
\Omega^\eps_n (\xi,\u):={}& 
\bma
\OM{t_{1,n}^\eps(\xi,\u,\cdot)}{t_{1,n}^\eps(\xi,\u,\cdot)} & \OM{t_{1,n}^\eps(\xi,\u,\cdot)}{t_{2,n}^\eps(\xi,\u,\cdot)}\\
\OM{t_{2,n}^\eps(\xi,\u,\cdot)}{t_{1,n}^\eps(\xi,\u,\cdot)} & \OM{t_{2,n}^\eps(\xi,\u,\cdot)}{t_{2,n}^\eps(\xi,\u,\cdot)}\\
\ema\nn\\
={}& \l({\g(\u)^3 m+ m^\eps_n(\xi,\u)}\r) \bma
 0 & 1\\
-1				& 0 
\ema.\nn
\ee 
Now we consider for any $(\bar\xi,\bar u)\in\R\times [-U- V(2),U+ V(2)]$, $(\bar v,\bar w)\in H^1(\R)\times L^2(\R)$ the matrix:
\be
{}&M^\eps_n(\bar\xi,\bar\u,\bar\v,\bar\w)\nn\\
={}&\bma
\Ltwo{\bma \dxi^2\p_n^\eps(\bar\xi,\bar\u,\cdot) \\ -\dxi^2\t_n^\eps(\bar\xi,\bar\u,\cdot) \ema}{\bma \bar\v(\cdot)\\ \bar\w(\cdot) \ema}
 & \Ltwo{\bma \du\dxi\p_n^\eps(\bar\xi,\bar\u,\cdot)\\ -\du\dxi\t_n^\eps(\bar\xi,\bar\u,\cdot) \ema}{\bma \bar\v(\cdot) \\ \bar w(\cdot) \ema} \\
\Ltwo{\bma \dxi\du\p_n^\eps(\bar\xi,\bar\u,\cdot)\\ -\dxi\du\t_n^\eps(\bar\xi,\bar\u,\cdot) \ema}{\bma \bar\v(\cdot) \\ \bar w(\cdot) \ema}  &  \Ltwo{\bma \du^2 \p_n^\eps(\bar\xi,\bar\u,\cdot)\\ -\du^2\t_n^\eps(\bar\xi,\bar\u,\cdot) \ema}{\bma \bar\v(\cdot) \\ \bar w(\cdot) \ema}
\ema.\nn
\ee
We use \cref{le controlvirtualtvproduct}, \cref{le transestimate} and H\"older's inequality similar to the proof of \cref{le uniform decomposition virtual} 
and obtain for all $(\bar\xi, \bar u)\in\R\times [-U- V(2),U+V(2)]$, $(\bar v,\bar w)\in H^1(\R)\times L^2(\R)$:
\be\la{M small v w}
 \Vert \l[\Omega^\eps_n (\bar\xi,\bar\u)\r]^{-1}  M^\eps_n(\bar\xi,\bar\u,\bar\v,\bar\w)\Vert \le C (\nhone{\bar\v}+ \nltwo{\bar\w}),
\ee
where we denote by $\Vert \cdot \Vert$ a matrix norm. Let $I=I_2$ be the identity matrix of dimension 2. Due to \re{M small v w} we are able to choose $\eps_0>0$ such that if $ \nhone{\bar v},\nltwo{\bar w}\le \eps_0$ then the matrix
$$
I+\l[\Omega^\eps_n (\bar\xi,\bar\u)\r]^{-1}  M^\eps_n(\bar\xi,\bar\u,\bar\v,\bar\w)
$$
is invertible by von Neumann's theorem. 
Using
\re{decomposition1 F virtual}-\re{decomposition2 F virtual} we express the orthogonality condition ${\cal N}^\eps(\t,\p,\xi,\u)=0$ from \cref{le uniform decomposition virtual} in terms of $(\v,\w,\xi,\u)$ 
and take its derivative with respect to $t$.
For simplicity of notation, we drop $(\t,\p,\xi,\u)$ and obtain
in matrix form:
\be
0{}&=\ddt\bma
{{\cal N}^\eps_1}\\
{{\cal N}^\eps_2}\\
\ema
=\Omega
\bma
\d\xi-u\\
\d\u - \luen (\xi,u)\\
\ema
+ M \bma
\d\xi-u\\
\d\u- \luen (\xi,u)\\
\ema
+ \bma
P_1\\
P_2
\ema,\nn
\ee
where 
$M=M^\eps_n(\xi,\u,\v,\w)$, $\Omega=\Omega^\eps_n (\xi,\u)$, $P_1=P_{1,n}^\eps(\xi,\u,\v,\w)$, $P_2 =P_{2,n}^\eps(\xi,\u,\v,\w)$,
\be
{}&
P_{1,n}^\eps(\xi,\u,\v,\w)\nn\\
=
{}&\int\dxi\ptn\w(x) -\dxi\ttn\Big(\dx^2\v(x) -\cos\ttn\v(x) \Big)\,dx\la{termsP11}\\
{}&+\int \u\dxi^2\ptn\v(x) -\u\dxi^2\ttn\w(x)  \,dx\la{termsP12}
%
%
%
%
\\
{}&+\int\du\dxi\ptn\v(x) -\du\dxi\ttn\w(x) \,dx\cdot \luen (\xi,u) 
{+\int 
\dxi\ptn
\Rtno
\,dx}
\nn\\
%
%
%
{}&{  -\int  
\dxi\ttn 
\Bigg(\fr{\sin\ttn \v^2(x)}{2}
+\tilde R(\v)(x)+\Rtnt\Bigg)
\,dx}\,\nn
\ee
and
\be
{}&P_{2,n}^\eps(\xi,\u,\v,\w)\nn\\
=
{}& \int \du\ptn  \w(x) -\du\ttn\Big(\dx^2\v(x)-\cos\ttn\v(x) \Big)\,dx\la{termsP21}\\ 
{}&+\int \u\dxi\du\ptn\v(x)-\u\dxi\du\ttn\w(x)\,dx\la{termsP22}\\
%
%
%
{}&{+
\int \du^2\ptn\v(x)-\du^2\ttn\w(x)\,dx \cdot \luen(\xi,u)}
{+\int
\du\ptn 
\Rtno
\,dx}\nn\\
{}&{-\int
\du\ttn 
\Bigg(\fr{\sin\ttn \v^2(x)}{2}+\tilde R(v)(x)+ \Rtnt \Bigg)\,dx}.\nn
\ee
If $ \nhone{v},\nltwo{ w}\le \eps_0$ then we obtain 
as mentioned above 
by von
Neumann's theorem that 
\be
\bma
\d\xi-u\\
\d\u- \luen (\xi,u)\\
\ema=
-\Big( I+ \Omega^{-1}M \Big)^{-1} [\Omega^{-1}P].\nn
\ee
We make a further assumption on $\eps_0$, namely that $\eps_0$ should be so small that the convergence rates in \cref{le: remainder uniformly bound} 
are satisfied.
Now we consider $P_1$ and $P_2$. 
The zeroth-order Taylor's approximations (in $\eps$) of expressions \re{termsP11}-\re{termsP12} and \re{termsP21}-\re{termsP22} respectively are  
$$
\Ltwo{\L\bma
\v(\cdot) \\ 
\w(\cdot) \\ 
\ema}{ 
\bma
-\dxi\p_0(\xi,\u,\cdot)\\
\dxi\t_0(\xi,\u,\cdot)\\
\ema}\text{ and }
\Ltwo{\L\bma
\v(\cdot) \\ 
\w(\cdot) \\ 
\ema}{ 
\bma
-\du\p_0(\xi,\u,\cdot)\\
\du\t_0(\xi,\u,\cdot)\\
\ema},$$
where $\L$
is given in \cref{de appendix Lxiu}.
Integration by parts and symplectic orthogonality yield that these approximations
vanish, which can also be deduced from \cref{le os onedim}.
Thus we obtain from \cref{le: remainder uniformly bound}, 
and similar arguments as above 
\be
|P_1|\le {}&C[\nhone{v}+\nltwo{w}] \eps^{k+1} + C\nhone{v}^2 +C \eps^{n+k+1} ,\nn\\
|P_2|\le {}&C[\nhone{v}+\nltwo{w}] \eps^{k+1} + C\nhone{v}^2 +C \eps^{n+k+1} .\nn
\ee
\epr
\section{Lyapunov Function}\la{ch:Lyapunov functional}
In this section we introduce the Lyapunov function and calculate its time derivative.
\bde
Let $(v,w)$ be given by \re{decomposition1 F virtual}-\re{decomposition2 F virtual}, with $(\xi,\u)$ obtained from \cref{le uniform decomposition virtual}.
We define the Lyapunov function $ L^\eps$ by 
\be 
 L^\eps {}&=\int \fr{\w^2(x)} 2 +\fr{(\dx\v(x))^2} 2 +\fr{\cos(\ttn) \v^2(x)}2+\u\w(x)\dx\v(x)\,dx\la{lyapunovfunction eps}\,
\ee
and the  auxiliary function $L$ by
\be 
 L {}&=\int \fr{\w^2(x)} 2 +\fr{(\dx\v(x))^2} 2 +\fr{\cos(\t_K(\Z)) \v^2(x)}2+\u\w(x)\dx\v(x)\,dx.\nn
\ee
\ede
\ble
\la{lelyapunovfunction}
It holds that
\be
\ddt L^\eps  
={}& (\u-\d\xi)\Bigg[\int 
-\u\dx\v(x)\Big\{-\dxi\ptn\Big\}
-\w(x)\Big\{-\dxi\ptn\Big\}\nn\\
{}&+[\cos(\ttn)\v(x)-\dx^2\v(x)]\dxi\ttn+\u\w(x)\dx\dxi\ttn 
\,dx\Bigg]\nn\\
{}&-(\d\u-\luen(\xi,u) )\Bigg[\int 
-\u\dx\v(x) \Big\{ -\du\ptn \Big\}-\w(x)\Big\{ -\du\ptn \Big\}\nn\\
{}&+[\cos(\ttn)\v(x)-\dx^2\v(x)]\du\ttn+\u\w(x)\dx\du\ttn \,dx\Bigg]\nn\\
{}&-\d\u\int \fr{\sin(\ttn)}2\du\ttn\v^2(x)\,dx
+(\d\xi -\u) \int  {\cos(\ttn)} \v(x)\dx\v(x)\,dx 
\nn\\
{}&
+\d\u\int \w(x)\dx\v(x)\,dx
+\int\w(x)[\fr{\sin\ttn\v^2(x)}{2}+\ti R(\v)(x) +\Rtnt ]\nn\\
{}&+\u\dx\v(x)[\fr{\sin\ttn\v^2(x)}{2}+\ti R(\v)(x) +\Rtnt  ] +\dx\v(x) \dx\Rtno \nn\\
{}& +\cos(\ttn)\v(x) \Rtno+u\w(x)\dx\Rtno\,dx.\nn
\ee
\ele
\bpr
We use a similar technique as in the proof of 
\cite[Lemma 2.1]{KoSpKu}.
We can assume that the initial data of our problem have compact support. This allows us to do the following computations (integration by parts etc.). The claim for non-compactly supported initial data follows by density arguments. We obtain the stated formula by taking the time derivative of \re{lyapunovfunction eps}, where we use  
\be
{}&\int\fr{\dt[\cos\ttn]}2\v^2(x)\,dx
\nn\\
={}&- \int \fr{\sin(\ttn)}2\d\xi\Big[\dxi\ttn+\dx\ttn\Big]\v^2(x)\,dx\nn\\
{}&+\int \d\xi {\cos(\ttn)} \v(x)\dx\v(x)-\d\u\fr{\sin(\ttn)}2\du\ttn\v^2(x)\,dx\nn
\ee
and 
$
\int \dx\v(x)\dx^2\v(x)+\w(x)\dx\w(x)\,dx=0 
$.
\epr
\section{Lower Bound}\la{ch: Lower bound}
Here we introduce a functional ${\cal  E}$ and prove a lower bound on ${\cal  E}$ by using symplectic orthogonality combined with functional analytic arguments. This will imply a lower bound on the Lyapunov function $L^\eps$,
which will play a key role in the proof of the main result.

\bde
For $(v,w)\in H^1(\R)\times L^2(\R)$,  $(\xi,u)\in\R\times (-1,1)$ we set
$$
{\cal  E}(\v,\w,\xi,\u):=
\fr 1 2 \int(\w(x)+\u\dx\v(x))^2+\v_Z^2(x)+\cos(\t_K(Z))\v^2(x) \,dx
,
$$
where $Z=\g(x-\xi)$ and $\v_Z(x)=\dZ\v(\fr Z \g +\xi)=\fr 1 {\g} \dx\v(x)$.
\ede
\noindent
A straightforward computation yields the following lemma.   
\ble \la{le Efunctional equal Lfunctional}
For $(v,w)\in H^1(\R)\times L^2(\R)$,  $(\xi,u)\in\R\times (-1,1)$ it holds that
\be
{\cal  E}(\v,\w,\xi,\u)=
\int \fr{\w^2(x)} 2 +\fr{(\dx\v(x))^2} 2 +\fr{\cos(\t_K(\Z)) \v^2(x)}2+\u\w(x)\dx\v(x)\,dx.\nn
\ee
\ele
\noindent
Recalling the relations \re{decomposition1 F virtual}-\re{decomposition2 F virtual}  we introduce a notation in order to be able to express the orthogonality conditions in terms of the variables $(\v,\w,\xi,\u)$ instead of the variables 
$(\t,\p,\xi,\u)$.
\bde
For $(v,w)\in H^1(\R)\times L^2(\R)$,  $(\xi,u)\in\R\times (-u_*,u_*)$ we set
\be 
{}&{\cal \check N}_1^\eps(\v,\w,\xi,\u)=\int \dxi\ptn \v(x) - \dxi\ttn\w(x)\,dx,\nn
\\
{}&{\cal \check N}_2^\eps(\v,\w,\xi,\u)=\int \du\ptn \v(x) - \du\ttn\w(x)\,dx.\nn
\ee
\ede
\noindent 
Now we prove a lower bound on the functional ${\cal  E}$. 
\ble \la{leStuart eps} 
Let $\eps_0>0$ be sufficiently small.
There exists $c>0$ such that if $\eps\in(0,\eps_0)$,
$(\xi,\u)\in [-\Xi,\Xi] \times [-U- V(2),U+ V(2)]\subset \R\times (-1,1)$ and $(v,w)\in H^1(\R)\times L^2(\R)$
satisfy
$$ {\cal \check N}_2^\eps(\v,\w,\xi,\u)=0$$ 
then
$$
{\cal  E}(\v,\w,\xi,\u)
\ge c (\nhone{v}^2 + \nltwo{ w}^2).
$$
\ele
\bpr We follow closely \cite{Stuart3} and \cite{Stuart1}. This proof is a slight modification of the proof of \cite[Lemma 4.3]{Stuart3}.
First of all we choose $\eps_0$ such that $\eps_0\in(0,\eps^*)$ with $\eps^*$ from \cref{def:virtual solitary manifold}. We will specify $\eps_0$ later.
Notice that the operator $-\dZ^2+\cos\t_K(Z)$ is nonnegative. It has (see \cite{Stuart2}) an one
dimensional null space spanned by $\t_K'(\cdot)$ and the essential spectrum $[1,\infty)$.
We argue by contradiction and assume that the result claimed is false: $ \forall j\in\N~\exists  \eps_j\in(0,\eps_0] ,~ (\xi_j,\u_j)\in [-\Xi,\Xi]\times [-U- V(2),U+ V(2)],
~
(\bar v_j,\bar w_j)\in H^1(\R)\times L^2(\R):$
\be  \la{notEnergy} 
\begin{split} 
{\cal \check  N}_2^{\eps_j}(\bar\v_j,\bar\w_j,\xi_j,\u_j)=0\,,~~~~
{\cal  E}(\bar\v_j,\bar\w_j,\xi,\u_j)
< \fr 1 j (\nhone{ \bar v_j}^2 + \nltwo{ \bar w_j}^2).
\end{split}
\ee 
This statement is also true for the sequences 
$ v_j:= {\bar v_j} {(\nhone{ \bar v_j}^2 + \nltwo{ \bar w_j}^2)^{-\fr 1 2}}$
and  
$w_j:=  {\bar w_j} {(\nhone{ \bar v_j}^2 + \nltwo{ \bar w_j}^2)^{-\fr 1 2}}$. 
Assuming that $\nltwo{  v_j}\OT{\to}{j\to\infty} 0$ we obtain
$\nltwo{( v_j)_x}\OT{\to}{j\to\infty} 0$ and 
$\nltwo{ w_j}\OT{\to}{j\to\infty} 0$. This is a contradiction to the fact that $\nhone{  v_j}^2 + \nltwo{  w_j}^2=1~\forall j\in\N$. By passing to a subsequence we may assume (without loss of generality) that there exists  ${\bar\delta}>0$ such that
\be\la{normvj}\nltwo{  v_{j}}^2\ge {\bar\delta} ~\forall j\in\N.\ee 
Since $(v_{j}, w_{j})$ is bounded in $H^1(\R)\times L^2(\R)$ we may assume that $ v_j \OT{\rightharpoonup}{H^1(\R)} v$ and $ w_j \OT{\rightharpoonup}{L^2(\R)} w$  by taking subsequences. 
Due to Rellich's theorem we may assume again by passing to subsequences that $ v_j \OT{\rightarrow}{L^2(\Omega)} v$, where $\Omega\subset\R$ is bounded and open. Passing to a further subsequence we assume almost everywhere convergence.
Due to the fact that
\be\la{rcosinus}
\exists~ r>0~~~\text{s.t.}~~~ |\cos(\t_K(Z))|>\fr 1 2 ~~~\text{for}~~~ |Z|>r
\ee 
and that $-\dZ^2+\cos\t_K(Z)$ is a nonnegative operator we obtain the estimate 
$$
{\cal  E}(\v_j,\w_j,\xi_j,\u_j)
\ge \fr {1} 4 \int_{-\infty}^{\fr {-r} {\g(u_j)}+\xi_j} \v_j^2(x) \,dx
+\fr {1} 4 \int_{\fr {r} {\g(u_j)}+\xi_j}^{\infty} \v_j^2(x) \,dx,
$$
where we used integration by parts and substitution.
We may extract a subsequence such that $\u_j \OT{\to}{\R}  u$, $\xi_j \OT{\to}{\R} \xi$ and $\eps_j \OT{\to}{\R}  \hat\eps$. 
It follows from \re{notEnergy} and from the previous estimate that 
$$
\int_{\{x\in\R:|x|\ge \ti r\}} \v_j^2(x) \,dx \OT{\to}{j\to\infty} 0
$$
for a sufficiently large $\ti r$. As a consequence, \re{normvj} and the strong convergence on bounded domains yield 
$\int_{\{x\in\R:|x| \le \ti r\}} v^2(x)\,dx\ge{\bar\delta},$ from which it follows that $v \not\equiv 0$. Weak convergence and the continuity of 
$
\eps \mapsto (\hat\t_n^\eps,\hat\p_n^\eps,\lambda_{n}^\eps)$
(see \cref{thimplicitfunctionIT1 alpha}) 
imply using the triangle inequality that
\be
 {\cal\check N}_2^{\hat\eps}(\v,\w,\xi,\u)={}&0\la{orth2eqzero}   
\ee
and
\be\la{Energyzero1}
\fr 1 2 \int\l(\w(x)+\u\v'(x)\r)^2\,dx
\le{}&\liminf_{j\to\infty}\fr 1 2 \int\l(\w_j(x)+\u_j\v_j'(x)\r)^2\,dx,\\
\fr 1 2 \int\l(\fr 1 {\g(u)}\v'(x)\r)^2\,dx
\le{}&\liminf_{j\to\infty}\fr 1 2 \int\l(\fr 1 {\g(u_j)}(\v_j)'(x)\r)^2\,dx.
\ee
Due to \re{rcosinus} we are able to apply Fatou's lemma for a sufficiently large $\ti r$ and obtain
\be
\la{Energyzero2}
\begin{split}
{}&\fr 1 2 \int_{\{x\in\R:|x| > \ti r\}} \cos(\t_K(\g(\u)(x-\xi)))\v^2(x)\,dx
\\
\le
{}&
\liminf_{j\to\infty}\fr 1 2 \int_{\{x\in\R:|x| > \ti r\}}\cos(\t_K(\g(\u_j)(x-\xi_j)))\v_j^2(x)\,dx,
\end{split}
\ee
where we have used that $(v_j)$ converges almost everywhere.
The dominated convergence theorem yields:
\be
\la{Energyzero3}\begin{split}
{}&
\fr 1 2 \int_{\{x\in\R:|x| \le \ti r\}} \cos(\t_K(\g(\u)(x-\xi)))\v^2(x)\,dx
\\
={}&\lim_{j\to\infty}\fr 1 2 \int_{\{x\in\R:|x| \le \ti r\}}\cos(\t_K(\g(\u_j)(x-\xi_j)))\v_j^2(x)\,dx.
\end{split}
\ee
\re{notEnergy} together with \re{Energyzero1}-\re{Energyzero3} imply that
$
{\cal  E}(\v,\w,\xi,\u)=0\,.
$
This yields $(v(x),w(x))=\eta(\t_K'(\g(u)(x-\xi)),-\u\g(u)\t_K''(\g(u)(x-\xi)))$ for some $\eta\not=0$, since $v \not\equiv 0$.
Using \cref{le transestimate} and the notation from \cref{decomposition of virtual manifold} we choose 
$\eps_0$ sufficiently small so that for all $ (\xi,\u)$
\be
\fr 1 {\g(u)}\l(\nltwo{\du\hat\t_n^{\hat\eps}(\xi,\u,x)} \nltwo{\dxi\p_0(\xi,\u,x)}
+\nltwo{\du\hat\p_n^{\hat\eps}(\xi,\u,x)} \nltwo{\dxi\t_0(\xi,\u,x)}\r)\le \fr m 2,\nn
\ee
which implies that ${\cal \check N}_2^{\hat\eps}(\v,\w,\xi,\u) \not= 0$. This yields a contradiction to \re{orth2eqzero}.
\epr
\bre
Let $(v,w)$ be given by \re{decomposition1 F virtual}-\re{decomposition2 F virtual}, with $(\xi,\u)$ obtained from \cref{le uniform decomposition virtual}. It holds that
$
L(t)={\cal  E}(\v(t),\w(t),\xi(t),\u(t)).
$
\ere
\section{Proof of \cref{maintheorem specialF kderivativesvanish virtual manifold}}\la{aprioriestimate F virtual manifold}
\subsection{Description of the Dynamics with Approximate
Equations for the Parameters $(\xi,\u)$}
We prove first a version of \cref{maintheorem specialF kderivativesvanish virtual manifold} with approximate equations for the parameters $(\xi,\u)$. 
\bth\la{th: prepairing main theorem}
Assume that the assumptions 
of \cref{maintheorem specialF kderivativesvanish virtual manifold}
on $n,k,\xi_s $ and $F$
are satisfied.
There exist $\eps_0,u_*,\tilde C >0$
and a map 
\be
(-\eps_0,\eps_0) \to Y_{2}^1(u_*),~
\eps \mapsto (\hat\t_n^\eps,\hat\p_n^\eps,\lambda_{ n}^\eps)\la{map in prepairing main theorem}
\ee 
of class $C^n$ such that the following holds. Let $\eps\in(0,\eps_0)$. 
Consider the Cauchy problem
\be\la{SGE1prepairing}
\partial_t \bma
\t \\
\p 
\ema
 =\l(\begin{matrix}
\p \\
\dx^2\t -\sin\t +\F\\
\end{matrix}\r)\,,
\bma
\t(0,x)\\
\p(0,x)
\ema =
\bma 
\t^\eps_n(\xi_s,\u_s,x)\\
\p^\eps_n(\xi_s,\u_s,x)
\ema 
+ \bma 
\v(0,x)\\
\w(0,x)
\ema,
\ee
where $(\t_n^\eps,\p_n^\eps,\lambda_{ n}^\eps)=(\t_0+\hatt_n^\eps ,\p_0+\hatp_n^\eps ,\lambda_{ n}^\eps)$ 
and $(\xi_s,\u_s)=(\xi(0),\u(0))\in \R\times (-1,1)$ such that the following assumptions are satisfied:
\begin{itemize}
\item[(a)] $|u_s|\le \tilde C\eps^{\fr{k+1}2}$.
\item[(b)] ${\cal N}^\eps(\t(0,\cdot),\p(0,\cdot),\xi_s,\u_s)=0$. 
\item[(c)] $\nhone{v(0, \cdot)}^2+\nltwo{w(0,\cdot )}^2\le \eps^{2n}$, where $(\v(0,\cdot),\w(0,\cdot))$ is given by \re{SGE1prepairing}.
\end{itemize}
Then the Cauchy problem defined by \re{SGE1prepairing}
has a unique solution on the time interval
\be
0 \le t \le T, ~\text{where}~ T=T(\eps,k )=\fr {1} {\tilde C\eps^{\beta(k )}}, ~~\beta(k )=\fr{k+1  }2.\nn
\ee
The solution may be written in the form
\be
\bma
\t(t,x)\\ 
\p(t,x)
\ema=\bma 
\t_n^\eps( \xi(t), \u(t),x)\\ 
\p_n^\eps( \xi(t), \u(t),x)
\ema+ 
\bma
\v(t,x)\\
\w(t,x)
\ema,\nn
\ee
where $\v, \w,\xi, \u$ have regularity
$
(v(t), w(t)) \in C^1([0,T ] , H^1(\R) \oplus L^2(\R)) 
$
and
$
(\xi(t), \u(t)) \in C^1([0,T ] ,\R \times (-1, 1)) 
$
such that the symplectic orthogonality condition
$$
{\cal N}^\eps(\t(t,\cdot),\p(t,\cdot),\xi(t),\u(t))=0\,
$$
is satisfied. There exist positive constants $c,C$ such that
\be
|\d\xi(t)-u(t)|  \le
C \eps^{n+k+1} ,\nn
~~~~~~	|\d\u(t)-\luen(\xi(t),\u(t))| \le{}&  C \eps^{n+k+1} ,\nn
\ee
and 
\be
\tnv{v}{0}{T}^2+\tnw{w}{0}{T}^2\le c \eps^{2n}.\nn
\ee
The constants $c,C$ depend on $F$ and $\xi_s$.
\eth
\noindent
Notice that the previous theorem describes the dynamics less precisely than \cref{maintheorem specialF kderivativesvanish virtual manifold}. However, in the previous theorem 
the orthogonality condition is satisfied which does not have to hold in \cref{maintheorem specialF kderivativesvanish virtual manifold}.\\
\noindent
The proof of \cref{th: prepairing main theorem} needs some preparation. The existence of $\eps_0>0$, $u_*>0$ and the map \re{map in prepairing main theorem} is ensured by \cref{thimplicitfunctionIT1 alpha}.
Now we suppose that \re{SGE1prepairing}
has a solution and we make some assumptions on $(v,w)$ given by \re{decomposition1 F virtual}-\re{decomposition2 F virtual} and on $(\xi,\u)$ obtained from \cref{le uniform decomposition virtual}. Then the  following lemma yields us more accurate information on $(v,w)$ and $(\xi,\u)$.  
\ble\la{le: prepairing main proof}
Assume that the assumptions 
of \cref{maintheorem specialF kderivativesvanish virtual manifold}
on $n,k,\xi_s,F$
are satisfied and let $0<\delta < 1/ {32}$. 
There exist $\eps_0, \tilde C>0$ such that the following statement holds. Let $\eps\in(0,\eps_0)$.
Assume that \re{SGE1prepairing}
has a  solution $(\t,\p) $ on $\l[0,\overline T\r]$ such that
$$
(\t,\p) \in C^1([0,\overline T], L^\infty(\R) \oplus L^2(\R)).
$$
Suppose that
$
0\le T \le t^*\le \overline T\,.
$
Suppose that $(v,w)$ is given by \re{decomposition1 F virtual}-\re{decomposition2 F virtual}, with $(\xi,\u)$ obtained from \cref{le uniform decomposition virtual} such that
$$
|u_s|\le \tilde C \eps^{\fr{k+1}2},~~~~ 
\tnv{v}{0}{T}^2+\tnw{w}{0}{T}^2\le  \eps^{2n-\delta}.$$
Then, provided
\be
0 \le T\le\fr {1} {\tilde C \eps^{ \beta(k )}}, ~~\beta(k )=\fr{k+1 }2,\nn
\ee
there exist $ c, C>0$ such that
\begin{itemize}
\item[(1)] $\forall  t \in [0,T]~(\xi(t),\u(t))\in\Sigma(5,\Xi),$
\item[(2)] $\tnv{v}{0}{T}^2+\tnw{w}{0}{T}^2\le \fr 1 {\CR c} (L(0)+C\eps^{2n}),$
where $c$ is from \cref{leStuart eps} and $C$ depends on $F,\xi_s$. 
\end{itemize}

\ele
\bre
Notice that the assumption $T \le t^*$ yields us the information:\\
$\forall  t \in [0,T]~ (\xi(t),\u(t))\in \Sigma(4,\Xi)\,.$
\ere
\bpr
Choose $\eps_0$ sufficiently small, in particular such that the lemmas used below can be applied and such that for any $\eps\in(0,\eps_0)$ the following statement holds: if $(v,w)\in H^1(\R)\times L^2(\R)$ satisfies $\nhone{v}^2+\nltwo{w}^2\le \eps^{2n-\delta}$ then it holds that $\nlinf{v}+\nltwo{w}\le \fr r 2$, where $r$ is from \cref{le uniform decomposition virtual}. This can be ensured by Morrey's embedding theorem.\\
\cref{le Modulation Equations} yields 
$\forall  t \in [0,T]$: 
\be
|\d\xi(t)-u(t)|  \le{}&  C[\nhone{v}+\nltwo{w}] \eps^{k+1} + C\nhone{v}^2 +C \eps^{n+k+1} 
\le
C \eps^{n+k+1-\delta},\nn
\\
	|\d\u(t)-\luen(\xi(t),\u(t))| \le{}&  C[\nhone{v}+\nltwo{w}] \eps^{k+1} + C\nhone{v}^2 +C \eps^{n+k+1} 
	\le
	C \eps^{n+k+1-\delta}.\nn
\ee
Then, using \cref{le: remainder uniformly bound}, it follows that there exists $\tilde C>1$ such that $\forall  t \in [0,T]$:
\be
 |u(t)-u(0)|
\le{}& \int _0^t |\d u(s)| \,ds 
\le
\tilde C \eps^{k+1} t \nn, 
\\
|\xi(t)-\xi(0)|
\le {}& \int _0^t |\d \xi(s)| \,ds 
%
\le 
\tilde C \eps^{n+k+1-\delta} t + \tilde C \eps^{k+1} t^2 + |u(0)|t.\nn
\ee
This implies (1) by choosing $\eps_0$ small enough and utilizing $|u_s|\le \tilde C \eps^{\fr{k+1}2}$.
By using \cref{leStuart eps}, \cref{le Efunctional equal Lfunctional}, \cref{lelyapunovfunction} and \cref{le: remainder uniformly bound} we obtain for times
$
0 \le t\le T\le  \fr 1 {\tilde C \eps^{\beta(k )}}, 
$
estimate (2):
\be
{}&
{\CR c}(\nhone{v(t)}^2+\nltwo{w(t)}^2)\nn
\\
\le 
{}&  
L(t)
=
L^\eps(t)+C\eps^{k+1} \tnv{v}{0}{t}^2
=
 L^\eps(0)
+\int_0^t \ddt{L^\eps}(t)\,dt +C\eps \tnv{v}{0}{t}^2\nn\\
 \le{}&  L^\eps(0)+C\eps^{2n}.\nn
\ee
\epr
\bth\la{th mainprooftrcontradiction kderivativesvanish}
Assume that the assumptions 
of \cref{maintheorem specialF kderivativesvanish virtual manifold}
on $n,k,\xi_s $ and $F$
are satisfied.
There exists $\eps_0,\tilde C >0$ such that the following statement holds. Let $\eps\in(0,\eps_0)$.
Assume that
\re{SGE1prepairing}
has a solution $(\t,\p)$  on $\l[0,\overline T\r]$ such that
$$
(\t,\p) \in C^1([0,\overline T], L^\infty(\R) \oplus L^2(\R)).
$$
Suppose that
$
0\le T \le \overline T\,
$
and that the assumptions (a),(b),(c)  of \cref{th: prepairing main theorem} are satisfied.  
Then, provided 
\be
0 \le T \le\fr {1} {\tilde C\eps^{\beta(k)}}, ~~\beta(k)=\fr{k+1 }2, \nn
\ee
it holds that $(v,w)$ given by \re{decomposition1 F virtual}-\re{decomposition2 F virtual} is well defined for times $[0,T]$ and 
there exists $\hat c>0$ such that
\begin{itemize}
\item[(1)]
$
\tnv{v}{0}{T}^2+\tnw{w}{0}{T}^2  \le \hat c \eps^{2n},
$
\item[(2)]
$
\forall  t \in [0,T]~(\xi(t),\u(t))\in \Sigma(5,\Xi).
$
\end{itemize}
\eth
\bpr
Let $\delta$ and $\tilde C$ be as in \cref{le: prepairing main proof}.
Choose $\eps_0$ sufficiently small, in particular such that
$$
\forall\eps\in(0,\eps_0):~\fr 2
{\CR c}( L(0)+C\eps^{2n}) < \eps^{2n-\delta},
$$
where $L(0)={\cal  E}(\v(0),\w(0),\xi_s,\u_s)$ and the constants $c,C$ are from \cref{le: prepairing main proof} (2).
Let $\eps\in(0,\eps_0)$.
Notice that $\Sigma(5,\Xi)
\subset \Sigma(4,\Xi)$.
We define an exit time
$$\ba
t_*:=\sup\Big\{{}& T>0:\tnv{v}{0}{t}^2+\tnw{w}{0}{t}^2\le \fr 2
{\CR c}( L(0)+C\eps^{2n}),\\ 
{}&  (\xi(t),\u(t))\in \Sigma(5,\Xi),~0\le t \le T\Big\}.
\ea $$
Suppose $t_* < \fr {1} {\tilde C\eps^{\beta(k )}}$. Then there exists a time $\hat t$ such that 
$
\fr {1} {\tilde C \eps^{\beta(k )}}>\hat t>t_*,
$ 
with
$$ \forall  t \in [0,\hat t]:~(\xi(t),\u(t))\in \Sigma(4,\Xi),~~~~(\xi(\hat t),\u(\hat t))\notin \Sigma(5,\Xi)$$ 
or
$$
\fr 1
{\CR c}( L(0)+C\eps^{2n})<\fr 2
{\CR c}( L(0)+C\eps^{2n})<\tnv{v}{0}{\hat t}^2+\tnw{w}{0}{\hat t}^2<\eps^{2n-\delta}.
$$
This leads a contradiction to \cref{le: prepairing main proof}. 
\epr
\noindent The previous theorem implies that the local solution of \re{SGE1prepairing}
discussed in \cref{local solution F virtual} is indeed continuable up to times $  {1} /({\tilde C\eps^{\beta(k)}})$ for $\eps\in (0,\eps_0)$. \cref{th mainprooftrcontradiction kderivativesvanish} and \cref{le Modulation Equations} yield the approximate equations for the parameters $(\xi,\u)$, which conclude the proof of \cref{th: prepairing main theorem}.
\subsection{ODE Analysis}
In this subsection we lay the groundwork for passing from the approximate equations for the parameters $(\xi,\u)$
in \cref{th: prepairing main theorem} to the ODEs given by \re{exactODE virtual1}.
We start with a preparing lemma.
\ble\la{le reference trajectory kderivativesvanish}
There exists $\eps_0>0$ such that the following statement holds. Let $\eps\in(0,\eps_0)$.
Let $\beta(k)=\fr{k+1}2$.
Let $\tixi=\tixi(s)$, $\tiu=\tiu(s)$, 
$\epsilon_1=\epsilon_1(s)$, $\epsilon_2=\epsilon_2(s)$ be $C^1$ real-valued 
functions. 
Suppose that   
$$|\epsilon_j(s)| \leq \bar c\eps^{n}$$ on $[0,T]$ for $j=1,2$.
Assume that on $[0,T]$, 
\be
\ds \tixi(s) ={}&  \tiu(s)+\epsilon_1(s),~~ \tixi(0)= \tixi_0,\nn\\
\ds \tiu(s) ={}&   \fr 1 {\eps^{2\beta(k )}} \lambda_{ n}^\eps(\tixi(s),\eps^{\beta(k )}\tiu(s)) +\epsilon_2(s),~~\tiu(0)=\tiu_0.\nn
\ee
Let $\hxi=\hxi(s)$ and $\hu=\hu(s)$ be $C^1$ real-valued 
functions which satisfy the exact equations
\be
\ds \hxi(s) ={}&  \hu(s) ,~~ \hxi(0)= \tixi_0,\nn\\
\ds \hu(s) ={}& \fr 1 {\eps^{2\beta(k )}} \lambda_{ n}^\eps(\hxi(s), \eps^{\beta(k )}\hu(s)),~~\hu(0)=\tiu_0.\nn
\ee
Then 
there exists $c>0$ such that the estimates
$$|\tixi(s)-\hxi(s)|\leq c\eps^{n}, 
\qquad |\tiu(s)-\hu(s)| \leq c\eps^{n},$$
hold on $[0,T]$.
\ele

\bpr
We follow very closely \cite[Lemma 6.1]{HoZwSolitonint}.  
We choose $\eps_0$ so small that the convergence rates in \cref{le: remainder uniformly bound} are satisfied for all $\eps\in(0,\eps_0)$. Let $\eps\in(0,\eps_0)$.
Let $x=x(s)$ and $y=y(s)$ be $C^1$ 
real-valued functions, $C\ge 1$, and let $(x,y)$ satisfy the differential 
inequalities:
$$
\left\{
\begin{aligned}
&|\dot x| \leq |y|, \\
&|\dot y| \leq C |x|+ C |y|,
\end{aligned}
\right. 
\qquad
\begin{aligned}
&x(0)=x_0,\\
&y(0)=y_0.
\end{aligned}
$$
For $z(s)=x^2+y^2$ the following estimate holds 
$$|\dot z| = |2x\dot x + 2y\dot y| \leq 2|x||y| 
+ 2C |x||y| +2C |y||y| \leq 4C(x^2+y^2) = 4Cz.$$ It follows from Gronwall's lemma that $z(s) \leq z(0)e^{4Cs}$.
Thus
\be  \begin{split} \la{Gronwall}
|x(s)| \leq \sqrt 2\max(|x_0|,|y_0|) \exp(2Cs),
~~~~|y(s)| \leq \sqrt 2\max(|x_0|,|y_0|) \exp(2Cs).
 \end{split}
\ee
Now we recall Duhamel's formula.
Let $X(s): \mathbb{R} \to \mathbb{R}^2$ be a two-vector function,  
$X_0\in \mathbb{R}^2$ a two-vector, and 
$A(s):\mathbb{R}\to (2\times 2\text{ matrices})$ a $2\times 2$ matrix function. We consider the ODE system $$\dot X(s) = A(s)X(s),~~~~X(s')=X_0$$ and denote its solution by $X(s)=S(s,s')X_0$ such that
\[ \frac{d}{ds} S(s,s')X_0 = A(s)S(s,s')X_0, \ \ S(s',s')X_0=X_0 . \]
Let $F(s):\mathbb{R}\to \mathbb{R}^2$ be a 2-vector function. We can express the 
solution of the inhomogeneous ODE system  
\begin{equation}
\dot X(s) = A(s)X(s) + F(s)\nn
\end{equation}
with initial condition $X(0)=0$ by Duhamel's formula
\begin{equation}
X(s) = \int_0^s S(s,s')F(s')ds'.\nn
\end{equation}
Let $U=\hu-\tiu$ and $\Xi=\hxi-\tixi$. These functions satisfy 
\be
\ds \Xi(s) =  U(s) + \epsilon_1(s), ~~~~~
\ds U(s) =  \fr 1 {\eps^{2\beta(k )}} 
\Big[\luen( \hxi(s),\eps^{\beta(k )}\hu(s))
-\luen( \tixi(s),\eps^{\beta(k )}\tiu(s))\Big] + \epsilon_2(s) .\nn
\ee
\noindent
Let 
\be
g(s)={}&\left\{
\begin{aligned}
& \fr 1 {\eps^{2\beta(k )}} \fr{\luen(\hxi(s),\eps^{\beta(k )}\hu(s))-\luen(\tixi(s),\eps^{\beta(k )}\hu(s))}{\hxi(s)-\tixi(s)}, & \ \ \text{if }\hxi(s) \neq \tixi(s),\\
& \ \ \fr 1{\eps^{2\beta(k )}} \partial_1\luen(\tixi(s),\eps^{\beta(k )}\hu(s)),  &\ \ \text{if } \hxi(s)=\tixi(s),
\end{aligned}\nn
\right.
%
\nn
\ee
\be
%
h(s)={}&\left\{
\begin{aligned}
&
\fr 1 {\eps^{2\beta(k )}} \fr{\luen(\tixi(s),\eps^{\beta(k )}\hu(s))-\luen(\tixi(s),\eps^{\beta(k )}\tiu(s))}{\tiu(s)-\hu(s)},
 & \ \ \text{if }\hu(s) \neq \tiu(s),\\
& \ \  \fr 1 {\eps^{\beta(k )}}\partial_2 \luen(\tixi(s),\eps^{\beta(k )}\hu(s)), 
&\ \ \text{if } \tiu(s)=\tiu(s).
\end{aligned}\nn
\right. 
\ee

\noindent We set
\be
A(s) =\bmat 0 & 1 \\ 
g(s) & h(s) 
\emat, 
\quad F(s) = \bmat \epsilon_1(s) \\ \epsilon_2(s) \emat, 
\quad X(s)=\bmat   \Xi(s) \\ U(s) \emat\,\nn
\ee
and obtain by Duhamel's formula:
\begin{equation}
\label{Duhamel}
\begin{bmatrix}
\Xi(s) \\ U(s)
\end{bmatrix}
= \int_0^s S(s,s') \begin{bmatrix} \epsilon_1(s') \\ \epsilon_2(s') 
\end{bmatrix} \, ds' .
\end{equation}
We use \cref{le: remainder uniformly bound} and apply \re{Gronwall} with
$$\begin{bmatrix} x(s) \\ y(s) \end{bmatrix} = S(s+s',s')\begin{bmatrix} 
\epsilon_1(s') \\ \epsilon_2(s') \end{bmatrix},  \quad \begin{bmatrix} x_0 
\\ y_0 \end{bmatrix} = \begin{bmatrix} \epsilon_1(s') \\ 
\epsilon_2(s') \end{bmatrix}.$$
It follows that 
$$\left| S(s,s') \begin{bmatrix} \epsilon_1(s') \\ 
\epsilon_2(s') \end{bmatrix} \right| \leq \sqrt 2 \begin{bmatrix} 
\exp(2C(s-s')) \\ \exp(2C(s-s')) \end{bmatrix} 
\max(|\epsilon_1(s')|,|\epsilon_2(s')|) .$$
Using \eqref{Duhamel} we obtain that on $[0,T]$
\be
&
|\Xi(s)| \leq \sqrt 2 \, T {\exp(2CT)}
\sup_{0\leq s\leq T}\max(|\epsilon_1(s)|,|\epsilon_2(s)|),\nn\\
&|U(s)| \leq \sqrt 2 \, T\exp( 2C T) 
\sup_{0\leq s\leq T}\max( |\epsilon_1(s)|,|\epsilon_2(s)|) ,\nn
\ee
which yields the claim.
\epr
\noindent 
In the following we show the relation between the parameters $(\xi,\u)$ selected by the implicit function theorem according to \cref{le uniform decomposition virtual} and the solutions $(\hat\xi,\hat\u)$ of the exact ODEs from the previous lemma.
\ble\la{le exact dynamics kderivativesvanish}
Assume that the assumptions 
of \cref{th: prepairing main theorem} are satisfied.
There exists $ \eps_0>0$ such that the following statement holds. Let $\eps\in(0, \eps_0)$,
$\beta(k )=\fr{k+1 }2$ and
$
s=\eps^{\beta(k )} t,
$
where 
$$
0\le s\le \fr 1 {\tilde C},~~~~0\le t\le \fr 1 {\tilde C\eps^{\beta(k )}}\,.
$$
Let $(\xi,\u)$ be the parameters selected according to \cref{le uniform decomposition virtual} and  
$(\hat\xi,\hat\u)$ from \cref{le reference trajectory kderivativesvanish}. 
Then there exists $c>0$ such that
\be
 {}&|\xi(t)-{\hxi(\eps^{\beta(k )} t)}|\le c\eps^{n }\,,
~~~~
|\u(t)-\eps^{\beta(k )}\hu(\eps^{\beta(k )} t)|\le c \eps^{n +\beta(k )}\,.\nn
\ee
\ele
\bpr
We choose $\eps_0$ as in \cref{th: prepairing main theorem}.
Let $\eps\in(0, \eps_0)$ and 
\be
\tixi(s)= \xi(s/\eps^{\beta(k )}),~~~~~ 
\tiu(s)= \fr 1 {\eps^{\beta(k )}} {u(s/\eps^{\beta(k )})}.\nn
\ee
For times
$
0\le t\le  ({{\tilde C}\eps^{\beta(k )}})^{-1}\,
$
\cref{th: prepairing main theorem} yields that 
\be
|\d\xi(t)-u(t)|  \le
C \eps^{n+k+1},~~~~~
%
%
%
	|\d\u(t)-\luen(\xi(t),\u(t))| \le
	C \eps^{n+k+1}.\nn
\ee
Thus $(\tixi,\tiu)$ satisfy the assumptions of \cref{le reference trajectory kderivativesvanish}, which implies that 
\be
{}&|\tixi(s)-\hxi(s)|=|\xi(t)-\hxi(\eps^{\beta(k )} t)|\le c\eps^{n },~~~~~
|\tiu(s)-\hu(s)|=|\fr{\u(t)}{\eps^{\beta(k )}}-\hu(\eps^{\beta(k )} t)|  \le c\eps^{n }.\nn
\ee 

\epr
\subsection{Completion of the Proof of \cref{maintheorem specialF kderivativesvanish virtual manifold}}

\cref{th: prepairing main theorem} yields the dynamics with the parameters $(\xi,\u)$  selected by the implicit function theorem according to \cref{le uniform decomposition virtual} on the time interval 
$
0\le t\le    ({\tilde C \eps^{ \beta(k )}})^{-1}\,.
$
Using \cref{le exact dynamics kderivativesvanish} and the the triangle inequality we can replace $(\xi(t),\u(t))$ 
with 
$
(\bar\xi(t),\bar\u(t)):=(\hxi(\eps^{\beta(k )} t),\eps^{\beta(k )}\hu(\eps^{\beta(k )} t))
$. The claim follows after possibly increasing the constant $\tilde C$ in the proof of \cref{th: prepairing main theorem}.
\qed

\appendix

\section{Preliminary Decompositions}\la{Preliminary Decompositions}
Let $\a,n\in\N$. Here we prove some decompositions for certain Sobolev spaces on $\R$ and on $\R^2$. 
We start with the spaces on $\R$ and prove an orthogonal decomposition of $H^{1}(\R)\oplus L^{2}(\R)$.
\bde\la{defoperatoren alpha}
We define the following spaces.
\begin{itemize}
\item [(a)] $L^{2,\alpha}(\R):=H^{0,\alpha}(\R) $.
\item [(b)] $H^{2,\a}_{\xi,\u,\perp}(\R):=\{\t\in H^{2,\a}(\R)~:~\ltwoa{\t(\cdot)}{\t_K'(\g(\cdot-\xi))}=0\}$.
\item [(c)] $L^{2,\a}_{\xi,\u,\perp}(\R):=\{\t\in L^2(\R)~:~\ltwoa{\t(\cdot)}{\t_K'(\g(\cdot-\xi))}=0\}$.
\end{itemize}
We define the following operators.
\begin{itemize}
\item [(d)] $L_{\xi,\u}^\a:  H^{2,\a}(\R)\subset L^{2,\a}(\R) \to L^{2,\a}(\R)$ given by\\
$$ (L_{\xi,\u}^\a \t)(x)= -(1-\u^2)\dx^2\t(x)+\cos(\t_K(\Z))\t(x).$$
\item [(e)] $M_{\xi,\u}^\a: H^{2,\a}(\R)\oplus \R\to L^{2,\a}(\R)$ given by \\
$$ \l(M_{\xi,\u}^\a\bma\t\\ \lambda \ema\r)(x) = (L_{\xi,\u} \t)(x)+\lambda \t_K'(\Z) .$$
\item [(f)] 
$\hat L_{\xi,\u}^\a= L_{\xi,\u}^\a\Big|_{H^{2,\a}_{\xi,\u,\perp}(\R)}$, $\hat M_{\xi,\u}^\a =  M_{\xi,\u}^\a\Big|_{H_{\xi,\u,\perp}^{2,\a}(\R)\oplus \R}$.
\end{itemize}
\ede

\ble\la{leranlaclosed}
$\ran L_{\xi,u}^\a$ is closed with respect to $L^{2,\alpha}(\R)$.
\ele
\bpr
We prove the case $\a=0$ that implies the claim. We consider the case $(\xi,u)=(0,0)$. The proof for a general $(\xi,u)\in \R\times(-1,1)$ works in the same way.
$L_{0,0}^0$ is self-adjoint and $0$ is an isolated point of $\sigma(L_{0,0}^0)$.
$l:=L_{0,0}^0|_{H^2(\R)\cap\langle\t_K'\rangle^\perp}$ is self-adjoint and has a bounded inverse (see \cite[Proposition 6.6]{HisSig}). Notice that $\ran L_{0,0}^0= \ran l$. Let $y_n=M x_n \OT{\to}{L^2} y $. Boundness yields  $x_n = l^{-1}y_n \OT{\to}{L^2} \overline{l^{-1}}y$, where $\overline{l^{-1}}$ denotes the bounded extension of $l^{-1}$ on the closure $\overline{\ran l}$. Since $l^*=l$ is a closed operator (see \cite[Proposition 4.9]{HisSig}), we obtain $l(\overline{l^{-1}}y)=y$.
\epr
\noindent
\cref{leranlaclosed} and the inverse mapping theorem imply the following lemma.
\ble
\begin{itemize} \la{le operator Lalpha R}
\item [(a)] $\ker L^\a_{\xi,\u} = \langle \t_K'(\g(\cdot-\xi)) \rangle$, $L^{2,\alpha}(\R)=\ran L^\a_{\xi,\u} \OT{\oplus}{L^{2,\a}} \langle \t_K'(\g(\cdot-\xi)) \rangle$.
\item [(b)] $L^{2,\alpha}(\R)=\ran \hat L^\a_{\xi,\u} \OT{\oplus}{L^{2,\a}} \langle \t_K'(\g(\cdot-\xi)) \rangle$.
\item [(c)] $\hat L^\a_{\xi,\u} \in L (H^{2,\alpha}_{\xi,\u,\perp}(\R),L^{2,\alpha}_{\xi,\u,\perp}(\R))$.
\item [(d)] $[\hat L^\a_{\xi,\u}]^{-1} \in L (L^{2,\alpha}_{\xi,\u,\perp}(\R),H^{2,\alpha}_{\xi,\u,\perp}(\R))$. 
\item [(e)] $\hat M^\a_{\xi,\u}\in L (H^{2,\alpha}_{\xi,\u,\perp}(\R)\oplus\R,L^{2,\alpha}(\R))$ and $\hat M^\a_{\xi,\u}$ is one-to-one and onto.
\end{itemize}
\ele
\bde\la{de appendix Lxiu}
We define the following operators.
\begin{itemize}
\item [(a)] $\L:  H^2(\R)\oplus H^1(\R)  \to H^1(\R) \oplus L^2(\R)$ given by\\
$$ \l(\L\tp\r) (x) =
\bma
-\u\dx\t(x)-\p(x) \\ 
-\dx^2\t(x)+\cos(\t_K(\Z))\t(x)-\u\dx\p(x) \\ 
\ema.
$$
\item [(b)] $\hatL:  \Big [ H^2(\R)\oplus H^1(\R) \Big ] \cap (\ker \L)^{\perp, L^2 \oplus L^2} \to H^1(\R)\oplus L^2(\R)$ given by\\ 
$$\l(\hatL\tp\r) (x) =
\bma
-\u\dx\t(x)-\p(x) \\ 
-\dx^2\t(x)+\cos(\t_K(\Z))\t(x)-\u\dx\p(x) \\ 
\ema.$$
\end{itemize}
\ede
\noindent 
\ble[orthogonal sum] \la{le os onedim}
\be
H^{1}(\R)\oplus L^{2}(\R) = \hatL\Bigg(\Big [ H^{2}(\R)\oplus H^{1}(\R) \Big ] \cap \langle t_\xi(\xi,\u,\cdot) \rangle^{\perp, L^2 \oplus L^2}\Bigg) \OT{\oplus}{\small L^2 \oplus L^2} \langle \J t_1(\xi,\u,\cdot) \rangle\,.\nn
\ee
\ele
\bpr
"$\supset$": clear.
"$\subset$": Let $( \bar v, \bar w)\in H^{1}(\R)\oplus L^{2}(\R)$. 
Orthogonal decomposition of $L^2(\R)\oplus L^2(\R)$ yields that there exists $\mu(\xi,\u)\in\R$ and $(\t_n, \p_n) \in H^{2}(\R)\oplus H^{1}(\R) $ such that
\be
\bma \bar v\\ \bar w\ema = { \scriptscriptstyle L^2 \oplus L^2- }\lim_{n\to\infty} \L\bma \t_n\\ \p_n\ema+ \mu(\xi,\u)  \J t_1(\xi,\u,\cdot)	,\nn
\ee
since
$\ker \L^* = \langle \J t_1(\xi,\u,\cdot) \rangle  
$ due to \cref{le operator Lalpha R}.
Hence $\bma  v\\ w\ema :={\scriptscriptstyle L^2 \oplus L^2 -}\lim_{n\to\infty} \L\bma \t_n\\ \p_n\ema\in\langle \J t_1(\xi,\u,\cdot) \rangle ^{\perp,L^2 \oplus L^2}\cap \Big[H^{1}(\R)\oplus L^{2}(\R)\Big]$ and thus
$
(-\u\g \v'(\cdot) +\g\w(\cdot)) \in \langle \t_K'(\g(\cdot-\xi) \rangle^{\perp , L^2}=\ran \hat L_{\xi,\u}
$
due to \cref{le operator Lalpha R} ($\a=0$).
By setting
$ 
\tit(x):= {[\hat L_{\xi,\u}]^{-1}}\Big(-\u\g \dx\v(x) +\g\w(x)\Big)
$,
$
\tip(x):=-\u\dx \tit(x)-\v(x)
$ we obtain that $(\ti\t,\ti \p) \in H^{2}(\R) \oplus H^{1}(\R)$ due to \cref{le operator Lalpha R} and that
$\hatL 
\bma 
\tit\\
\tip
\ema_{(\ker \L)^{\perp, L^2 \oplus L^2}}
 = \bma
\v\\
\w
\ema$.
\epr
\noindent
We turn to decompositions of Sobolev spaces on $\R^2$. Here
\cref{le os onedim} will be needed later in the proofs. We start with a definition and some preparing lemmas.
Unlike the one-dimensional case we consider $\xi$ not as
a fixed parameter anymore, but as a new variable.
\bde\la{defoperatoren alpha}
We define the following spaces.
\begin{itemize}
\item [(a)] $H^{2,\a}_{\perp}(\R^2):=\{\t\in H^{2,\a}(\R^2)~|~\forall\lambda\in H^2(\R):~\ltwortwoaxiZ{\t(\xi,Z)}{\lambda(\xi)\t_K'(Z)}=0\}$.
\item [(b)] $H^{2,\a}_{\u,\perp}(\R^2):=\{\t\in H^{2,\a}(\R^2)~|~\forall\lambda\in H^{2,\a}(\R):~\ltwortwoaxix{\t(\xi,x)}{\lambda(\xi)\t_K'(\g(x-\xi))}=0\}$.
\item [(c)] $L^{2,\a}_{\u,\perp}(\R^2):=\{\t\in L^{2,\a}(\R^2)~|~\forall\lambda\in H^{2,\a}(\R):~\ltwortwoaxix{\t(\xi,x)}{\lambda(\xi)\t_K'(\g(x-\xi))}=0\}$.
\item [(d)] $\bar {\cal Y}^\a$ is the space $H^{2,\a}(\R^2)\oplus H^{1,\a}(\R^2)$
with the finite norm
$$|y|_{\bar {\cal Y}^\a} = |\t|_{H^{2,\a}(\R^2)}+ |\lambda |_{H^{1,\a}(\R^2)}. $$
\item [(e)] $\bar {\cal Z}^\a$ is the space $L^{2,\a}(\R^2)$ 
with the finite norm
$$|z|_{\bar {\cal Z}^\a} = |z|_{L^{2,\a}(\R^2)}.$$
\item [(f)] ${\cal Y}^\a={\cal Y}^\a(\u_*)$ is the space \\
$\ba {}& \bigg\{ y=(\t,\lambda ) \in C(I(\u_*), \bar {\cal Y}^\a) : \Vert y \Vert_{{\cal Y}^\a(\u_*)} <\infty,~\forall \u\in I(\u_*):
\t(\u)\in H^{2,\a}_{\u,\perp}(\R^2) \bigg\}\ea\\$\\
with the finite norm
$$\Vert y \Vert_{{\cal Y}^\a(\u_*)} =\sup_{\u\in I(\u_*)}  |y|_{\bar {\cal Y}^\a}\,.$$
\item [(g)] ${\cal Z}^\a={\cal Z}^\a(\u_*)$ is the space 
$\bigg\{ z \in C( I(\u_*), \bar {\cal Z}^\a) : \Vert z \Vert_{{\cal Z}^\a(\u_*)} <\infty \bigg\}\,$ 
with the finite norm
$$
\Vert z \Vert_{{\cal Z}^\a(\u_*)} =\sup_{\u\in I(\u_*)}  |z|_{\bar {\cal Z}^\a}.
$$
\end{itemize}
We define the following operators.
\begin{itemize}
\item [(h)] $L^\a: H^{2,\a}(\R^2)\subset L^{2,\a}(\R^2) \to L^{2,\a}(\R^2)$ given by\\ $$(L^\a \t)(\xi,Z) = -\dZ^2\t(\xi,Z)+\cos(\t_K(Z))\t(\xi,Z).$$
\item [(i)] $L_{\u}^\a:  H^{2,\a}(\R^2)\subset L^{2,\a}(\R^2) \to L^{2,\a}(\R^2)$ given by\\
$$(L_{\u}^\a \t)(\xi,x) = -(1-\u^2)\dx^2\t(\xi,x)+\cos(\t_K(\Z))\t(\xi,x).$$
\item [(j)] $M_{\u}^\a: H^{2,\a}(\R^2)\oplus H^{2,\a}(\R)\to L^{2,\a}(\R)$ given by \\
$$\l(M_{\u}^\a \bma\t \\ \lambda  \ema\r) (\xi,x) = ( L_{\u}^\a \t)(\xi,x)+\lambda(\xi) \t_K'(\Z) .$$
\item [(k)] $M^\a: H^{2,\a}(\R^2)\oplus H^{2,\a}(\R)\to L^{2,\a}(\R)$ given by\\
$$\l(M^\a \bma\t\\ \lambda \ema\r)(\xi,Z) = (L^\a \t)(\xi,Z)+\lambda(\xi) \t_K'(Z) .$$
\item [(l)] $
\hat L^\a=L^\a\Bigg|_{ H^{2,\a}_{\perp}(\R^2)} 
$, $\hat L_{\u}^\a=L_{\u}^\a\Bigg|_{H^{2,\a}_{\u,\perp}(\R^2)}$,
$\hat M_{\u}^\a=M_{\u}^\a\Bigg|_{H^{2,\a}_{\u,\perp}(\R^2)\oplus H^{2,\a}(\R)}$,
$\hat M^\a=M^\a\Bigg|_{H^{2,\a}_{\perp}(\R^2)\oplus H^{2,\a}(\R)}$.
\end{itemize}
\ede

\ble\la{le kernelEV alpha}
\begin{itemize} 
\item [(a)] $\ker L_u^\a= \{\ \t\in H^{2,\a}(\R^2)  ~|~	\t(\xi,x)=\lambda(\xi)\t_K'(\g(x-\xi)),~\lambda\in H^{2,\a}(\R)\}$.
\item [(b)] $0$ is an isolated eigenvalue of $L_u^\a$.
\end{itemize}
\ele
\bpr
The claim follows from \cref{le operator Lalpha R} combined with the fact that the operator $-\dZ^2+\cos\t_K(Z)$ is nonnegative with essential spectrum $[1,\infty)$ and with one
dimensional null space spanned by $\t_K'(\cdot)$ (see \cite{Stuart2}).
\epr
\noindent
Using that 
$L_u$ is self-adjoint and $0$ is an isolated point of $\sigma(L_u)$ one proves 
analogously to \cref{leranlaclosed}
the following lemma.
\ble \la{le: ranLalpha closed}
$\ran L_u^\a$ is closed with respect to $L^{2,\alpha}(\R^2)$.
\ele
\noindent
Analogously to 
\cref{le operator Lalpha R} we obtain the next lemma.
\ble\la{le operator LRtwo alpha}
\begin{itemize} 
\item [(a)] $L^{2,\alpha}(\R^2)=\ran L^\a_{\u} \OT{\oplus}{L^{2,\a}} \ker L_u^\a $.
\item [(b)] $L^{2,\alpha}(\R^2)=\ran \hat L^\a_{\u} \OT{\oplus}{L^{2,\a}} \ker L_u^\a $.
\item [(c)] $\hat L^\a_{\u} \in L (H^{2,\alpha}_{\u,\perp}(\R^2),L^{2,\alpha}_{\u,\perp}(\R^2))$.
\item [(d)] $\l[\hat L^\a_{\u}\r]^{-1} \in L (L^{2,\alpha}_{\u,\perp}(\R^2),H^{2,\alpha}_{\u,\perp}(\R^2))$. 
\item [(e)] $\hat M_{\u}^\a,\hat M^\a$ are one-to-one, onto, bounded and the inverse mappings are also bounded.
\end{itemize}
\ele

\ble\la{le smallOPRtwo alpha}
Let ${\frak m}^\a:{\cal Y}^\a \to {\cal Z}^\a,~
(\t,\lambda) \mapsto {\frak m}^\a(\t,\lambda)$ be the linear operator, given by
\be  {\frak m}^\a(\t,\lambda)(\u) =\hat M_{\u}^\a (\t(\u),\lambda(\u)).\nn
\ee
Then ${\frak m}^\a$ is one-to-one, onto and bounded, i.e., $\l[{\frak m}^\a\r]^{-1}$ is bounded.
\ele
\bpr 
${\frak m}^\a$ is well defined and ${\frak m}^\a$
one-to-one due to
\cref{le operator LRtwo alpha}. In order to see that 
${\frak m}^\a$ is onto let $v \in {\cal Z}^\a$. Due to 
\cref{le operator LRtwo alpha} for each $\u\in I$ there exists $(\t(\u) ,\lambda(\u))\in H_{\u,\perp}^{2,\a}(\R^2)\oplus H^{2,\a}(\R) $ such that 
\be\la{transformedMalpha}\
\v(u)(\xi,\fr Z \g +\xi)=\hat M^\a  \bma \bar\t(u) \\ \lambda(u) \ema (\xi,Z) =\l(\hat L^\a   \bar\t(u)\r)(\xi, Z )  +\lambda(u)(\xi)  \t_K'(Z),
\ee
where $\bar\t(u)(\xi,Z)=\t(u)(\xi,\fr Z \g +\xi)$. 
It holds for $h\in H^2(\R^2)$ the inequality
\be
{}&\nhtwortwoxix{h(\xi,x)}\le 
\sqrt  5\g(u)^{ \fr 32}\nhtwortwoxiZ{h(\xi,\fr Z {\g(u)} +\xi)}\la{contprep}.
\ee
Using \re{contprep} for $h(\xi,x)=\l(1+|\xi|^2+|x|^2\r)^\fr\a 2\l[\t(u)-\t(\bar\u)\r](\xi,x)$, \re{transformedMalpha} and \cref{le operator LRtwo alpha} we obtain:
\be
{}&\l| \bma \t(u)-\t(\bar u) \\ \lambda(u)-\lambda(\bar u) \ema\r|_{\bar {\cal Y}^\a}\nn\\
\le {}& \g(u)^{ \fr 32}C(\a) \l\Vert \l[\hat M^\a\r]^{-1} \r\Vert
\nltwortwoxiZ{\l(1+|\xi|^2+| Z |^2 \r)^{\fr\a 2}\l[\v(\u)(\xi,\fr Z {\g(u)} +\xi)
-\v(\u)(\xi,\fr Z {\g(\bar\u)} +\xi)\r]}\nn\\
{}&+
 \g(u)^{ \fr 32} C(\a) \l\Vert \l[\hat M^\a\r]^{-1} \r\Vert
\nltwortwoxiZ{\l(1+|\xi|^2+| Z |^2 \r)^{\fr\a 2}\l[
\v(\u)(\xi,\fr Z {\g(\bar\u)} +\xi)
-\v(\bar\u)(\xi,\fr Z {\g(\bar\u)} +\xi)\r]}\nn\\
{}&+\g(u)^{ \fr 32}C(\a) \nhtwortwoxiZ{
\l(1+|\xi|^2+| Z |^2 \r)^{\fr\a 2}\l[
\t(\bar u)(\xi,\fr Z {\g(\bar u)} +\xi)-\t(\bar\u)(\xi,\fr Z {\g(u)} +\xi)
\r]} .\nn
\ee
This implies that
$
(\t,\lambda)\in {\cal Y}^\a
$,
since $v \in {\cal Z}^\a$. 
The inverse mapping theorem yields that $\l[{\frak m}^\a\r]^{-1}$ is bounded, since ${\frak m}^\a$ is bounded.
\epr
\noindent In the following we introduce the operator $\hat {\cal L}_u^\a$ that will be used for the main decomposition of this appendix in \cref{co ds Rtwo alpha}.  
\bde We define the following operators.
\begin{itemize}
\item [(a)] ${\cal L}_u^\a:  H^{2,\a}(\R^2)\oplus H^{1,\a}(\R^2)  \to H^{1,\a}(\R^2)\oplus L^{2,\a}(\R^2)$ given by\\
$$
\l({\cal L}_u^\a
\tp \r)(\xi,x) =
\bma
-\u\dx\t(\xi,x) -\p(\xi,x)  \\ 
-\dx^2\t(\xi,x) +\cos(\t_K(\Z))\t(\xi,x) -\u\dx\p(\xi,x)  \\ 
\ema
.$$
\item [(b)] $\hat{\cal L}_u^\a=\hat{\cal L}_u^\a\Bigg|_{  \Big [ H^{2,\a}(\R^2)\oplus H^{1,\a}(\R^2) \Big ] \cap \l(\ker {\cal L}_u\r)^{\perp, L^{2,\a} \oplus L^{2,\a}}} $.
\end{itemize}
\ede

\ble[orthogonal sum] \la{le os Rtwo alpha}
\be{}&H^{1,\a}(\R^2)\oplus L^{2,\a}(\R^2)\nn\\
 = {}&\hat{\cal L}_u^\a\Bigg(\Big [ H^{2,\a}(\R^2)\oplus H^{1,\a}(\R^2) \Big ] \cap \l(\ker {\cal L}_u^\a\r)^{\perp, L^{2,\a} \oplus L^{2,\a} }\Bigg) \OT{\oplus}{\small L^{2,\a}  \oplus L^{2,\a} } \{ \lambda  \J t_1 (\u),\lambda\in H^{2,\a}(\R)\}.\nn
\ee
\ele
\bpr 
It follows from \cref{le kernelEV alpha} that $\ker \l[{\cal L}_u^\a \r]^*
=\{ \lambda  \J t_1 (\u),\lambda\in H^{2,\a}(\R)\}$. One proves first the case $\a=0$ analogously to the proof of \cref{le os onedim} which can be used to deduce the case $\a\not=0$.
\epr
\bco[direct sum]\la{co ds Rtwo alpha}

\be
{}& H^{1,\a}(\R^2)\oplus L^{2,\a}(\R^2)\nn\\
={}& \hat{\cal L}_u^\a \Bigg(\Big [ H^{2,\a}(\R^2)\oplus H^{1,\a}(\R^2) \Big ] \cap \l(\ker {\cal L}_u^\a\r)^{\perp, L^{2,\a} \oplus L^{2,\a}}\Bigg) \OT{\oplus}{\small L^{2,\a}  \oplus L^{2,\a} } \{ \lambda t_2 (\u),\lambda\in H^{2,\a}_\xi(\R)\}\,.\nn
\ee
\eco  
\bpr
"$\supset$": clear. 
"$\subset$": Let $(v,w)\in H^{1,\a}(\R^2)\oplus L^{2,\a}(\R^2)$ then
there exists due to \cref{le os Rtwo alpha} $( \t, \p) = ( \t(u), \p(u)) \in H^{2,\alpha}(\R^2)  \oplus H^{1,\alpha}(\R^2) \cap \l(\ker {\cal L}_u^\a\r)^{\perp, L^{2,\a} \oplus L^{2,\a}}$ and $\lambda=\lambda(u)\in H^{2,\a}(\R)$ such that 
$$
\bma
\v\\
\w
\ema
=\hat{\cal L}_u^\a \bma\t\\ \p\ema+\lambda \J t_1(u).
$$ 
%
Assume without loss of generality 
$\nhtwoa\lambda\not=0$, then
$
\Ltwortwoa{\lambda  t_1 (u)}{\lambda  \J t_2 (u)}
\not= 0\,.
$
Thus
due to \cref{le os Rtwo alpha} there exist
$( \bar\t , \bar\p ) =( \bar\t(\u) , \bar\p(\u)) \in H^{2,\alpha}(\R^2)  \oplus H^{1,\alpha}(\R^2) \cap \l(\ker {\cal L}_u^\a\r)^{\perp, L^{2,\a} \oplus L^{2,\a}}$ and 
$0\not=\bar\lambda=\bar\lambda(\u)\in H^{2,\a}(\R)$ such that
\be
\la{directsum alpha}
\lambda(\u) t_2(\u)={\cal L}_u^\a  \bma\bar\t(\u)\\ \bar\p(\u)\ema+\bar\lambda(\u) \J t_1(\u).
\ee
This is an identity in $H^{1,\a}(\R^2)\oplus L^{2,\a}(\R^2)$. 
Fixing $\xi$ and pairing this identity with $\J t_1(\xi,\u,\cdot)$ in
$L_x^2(\R)\oplus L_x^2(\R)$ yields due to \cref{le os onedim} for a.e. $\xi\in \R$
the identity
$  
\lambda(\xi,\u)=	\eta(u) \bar\lambda(\xi,\u)\,,
$ 
where 
$ \eta(u):={\g(u)^{-3}m^{-1}}\l({\u^2\g^3\nltwo{\t_K''}^2+
\g\nltwo{\t_K'}^2}\r)\in \R\,.
$ 
Thus using \re{directsum alpha} we obtain
\be
\bma
\v \\
\w 
\ema
=\hat{\cal L}_u^\a \bma\t(\u)\\ \p(\u)\ema+\lambda(\u) \J t_1(\u)
=\hat{\cal L}_u^\a\l( \bma\t(\u)\\ \p(\u)\ema- {\eta(u)}\bma\bar\t(\u)\\ \bar\p(\u)\ema\r)+ {\eta(u)} \lambda(\u) t_2(\u).\nn
\ee
The sum is direct
due to \re{directsum alpha}.
\epr
\section{Proof of \cref{le invertibilityMxiCtwo alpha}}\la{Invertibility of the Linearization}
Let $\a,n\in\N$. We want to show that the linear operator $
{\frak M}_n^\a:  Y_n^\a(u_*)  \to Z_n^\a(u_*)
$
is invertible if $u_*$ is small. The operator ${\frak M}_n^\a$ contains derivatives with respect to $\xi$ and $x$ which makes it difficult to analyze it. Therefore we consider first an operator $\ti{\frak M}^\a:  Y^\a \to Z^\a,~
$
which contains only derivatives with respect to $x$. Using \cref{co ds Rtwo alpha} we prove that ${\frak M}_n^\a$ is invertible.

\bde We define the following operators.
\begin{itemize}
\item [(a)] ${\cal M}_u^\a:  H^{2,\a}(\R^2)\oplus H^{1,\a}(\R^2) \oplus H^{2,\a}(\R)\to H^{1,\a}(\R^2)\oplus L^{2,\a}(\R^2) $ given by\\
$$\l( {\cal M}_u^\a ( \t,
\p,
\lambda ) \r)(\xi,x) =
\l({\cal L}_{\u}^\a\tp\r) (\xi,x)+ \lambda(\xi) 
t_2(\xi,\u,x).$$
\item [(b)] 
$ {\cal \hat M}_u^\a= {\cal M}_u^\a\Bigg|_{\Big [ H^{2,\a}(\R^2)\oplus H^{1,\a}(\R^2) \Big ] \cap \l(\ker {\cal L}_u^\a\r)^{\perp, L^{2,\a} \oplus L^{2,\a}} \oplus H^{2,\a}(\R) } $.\\
\end{itemize}
We define the following spaces.
\begin{itemize}
\item [(c)] $\bar Y^\a$ is the space $H^{2,\a}(\R^2) \oplus H^{1,\a}(\R^2) \oplus H^{2,\a}(\R)$
with the finite norm
$$
|y|_{\bar Y^\a} = |\t|_{H^{2,\a}(\R^2)}+ |\p|_{H^{1,\a}(\R^2)}+|\lambda |_{H^{2,\a}(\R)}.
$$ 
\item [(d)] $\bar Z^\a$ is the space $H^{1,\a}(\R^2) \oplus L^{2,\a}(\R^2)$ 
with the finite norm
$$
|z|_{\bar Z^\a} = |\v|_{H^{1,\a}(\R^2)}+ |\w|_{L^{2,\a}(\R^2)}.
$$
\item [(e)] $ Y^\a  = Y^\a(\u_*)$ is the space\\
$\ba {}& \bigg\{  y=(\t,\p,\lambda) \in C( I(\u_*), \bar Y^\a) : \Vert y \Vert_{Y^\a(\u_*)} <\infty;~\forall~ \u\in I(\u_*),~\forall~\mu\in H^{2,\a}(\R):\\
{}& \Ltwortwoaxix{ \bma \t(\u)(\xi,x)\\
\p(\u)(\xi,x) \ema}{\mu(\xi)\bma \t_K'(\Z)\\
-\u\g\t_K''(\Z)\ema} =0 \bigg\}\,\ea\\$\\
with the finite norm
$$\Vert y \Vert_{Y^\a(\u_*)} =\sup_{\u\in I(\u_*)}  |y|_{\bar Y^\a}.$$
\item [(f)] $Z^\a=Z^\a(\u_*)$ is the space $\bigg\{ z =(\v,\w) \in C(I(\u_*), \bar Z^\a) : \Vert z \Vert_{Z^\a(\u_*)} <\infty \bigg\}\,$ 
with the finite norm
$$
\Vert z \Vert_{Z^\a(\u_*)} =\sup_{\u\in I(\u_*)}  |z|_{\bar Z^\a}.
$$
\end{itemize}
\ede
\ble\la{le invertibilityMx alpha}
The linear operator
$
\ti{\frak M}^\a:  Y^\a(\u_*) \to Z^\a(\u_*),~
( \t,
\p,
\lambda ) \mapsto \ti{\frak M}^\a ( \t,
\p,
\lambda )\,,
$
given by
\be
{}&\ti{\frak M}^\a
(\t,
\p,
\lambda )(\u)=
 {\cal \hat M}_u^\a ( \t(\u),
\p(\u),
\lambda  (\u)),\nn
\ee
is invertible.
\ele
\bpr
$\ti{\frak M}^\a$ is onto: Let $(
\v,
\w
)\in Z^\a$. Due to \cref{co ds Rtwo alpha} for all $\u\in I$ there exist 
$
( \t(\u),\p(\u)) \in [ H^{2,\a}(\R^2)\oplus H^{1,\a}(\R^2) \Big ] \cap \l(\ker {\cal L}_u^\a\r)^{\perp, L^{2,\a} \oplus L^{2,\a}}
$ 
and $\lambda (\u)\in H^{2,\a}(\R)$ such that
\be \la{operatorLalphaInTheProof}
\bma
\v(\u)\\
\w(\u)
\ema
={\cal \hat L}_u^\a\bma \t(\u)\\ \p(\u) \ema
+\lambda(\u)t_2(u).
\ee
This is an identity in $H^{1,\a}(\R^2)\oplus L^{2,\a}(\R^2)$. 
By fixing $\xi$, using \cref{le os onedim} and pairing \re{operatorLalphaInTheProof} with $ \J t_\xi(\xi,\u,x)$
in $L_x^{2,\a}(\R)\oplus L_x^{2,\a}(\R)$ we obtain 
$\lambda\in  C (I , L^{2,\a}(\R))$.
Further \re{operatorLalphaInTheProof} yields 
\be
\bma
 ((1+\u^2)^{-1}\g^{-3}\t )_{\ker L_{\u}^\a} \\
\lambda \\
\p
\ema
{}&=\bma
\l[\hat M_u^\a\r]^{-1}\Big[  {(1+\u^2)^{-1}\g^{-3}} 
\l(\w-\u\dx\v- 2\lambda[u^2\g^4(x-\xi)\t_K''(Z)]\r) \Big]\\
-\u\dx\t-\v+\lambda[\u\g^3(x-\xi)\t_K'(Z)]
\ema,\nn
\ee
where $Z=\Z$.
Hence $ ( ((1+\u^2)^{-1}\g^{-3}\t )_{\ker L_{\u}^\a} ,
\lambda  )\in {\cal Y}^\a$ due to \cref{le smallOPRtwo alpha}, since $  (1+\u^2)^{-1}\g^{-3}\l(\w-\u\dx\v\r) \in {\cal Z}^\a$ and $\lambda\in  C (I , L^{2,\a}(\R))$.
Thus
$
(\t,\p,\lambda) \in Y^\a
$ and $\ti{\frak M}^\a(\t,\p,\lambda)=(\v,\w) $.
$\ti{\frak M}^\a$ is one-to-one due to \cref {co ds Rtwo alpha}. The inverse mapping theorem yields the claim.
\epr
\noindent Next, we want to show that the operator norm of 
$ {\cal \hat M}_u^{-1} $
is bounded by a function, which is continuous  in $u$. We start with a preparing lemma.
\ble[Norm of ${\l[\hat M_u^\a\r]}^{-1}$] There exists a constant $c^\a>0$ such that
\la{le Muinvestimate alpha}
\be
\l\Vert \l[\hat M_u^\a\r]^{-1} \r\Vert_{L(L^{2,\a}(\R^2),H^{2,\a}(\R^2)\oplus H^{2,\a}(\R) )}\le  \g(u) c^\a \l\Vert \l[ \hat M^\a\r]^{-1} \r\Vert_{L(L^{2,\a}(\R^2),H^{2,\a}(\R^2)\oplus H^{2,\a}(\R) )}.\nn
\ee
\ele
\bpr
Let
$
\nltwortwoa{v}\le 1.
$
Due to \cref{le operator LRtwo alpha} and \cref{le kernelEV alpha} there exists  
$
(\t , \lambda ) \in H_{\u,\perp}^{2,\a}(\R^2)\oplus H^{2,\a}(\R), 
$
such that
\be 
\v(u)(\xi,\fr Z \g +\xi)=\hat M^\a  \bma \bar\t(u) \\ \lambda(u) \ema (\xi,Z) =\l(\hat L^\a   \bar\t(u)\r)(\xi, Z )  +\lambda(u)(\xi)  \t_K'(Z),\nn
\ee
where $\bar\t(u)(\xi,Z)=\t(u)(\xi,\fr Z \g +\xi)$. 
Using \re{contprep}
for $h(\xi,x)=\l(1+|\xi|^2+|x|^2\r)^\fr\a 2 \t(\xi,x)$
we obtain 
\be
{}&\l|\l(\l[\hat M_u^\a\r]^{-1}\v\r)(\xi,x)\r|_{H^{2,\a}_{\xi,x}(\R^2)\oplus H^{2,\a}_{\xi}(\R)}
\le 
\sqrt 5\g(u) \l\Vert\l[ \hat M^\a\r]^{-1} \r\Vert
\nltwortwoaxiZ{\v(\xi,Z)}.\nn
\ee
\epr
\ble[Norm of $ {\l[{\cal \hat M}_u^\a\r]}^{-1} $]
\la{le calMuinvestimate alpha}
There exists a continuous function $C^\a : (-1,1)\to \R$ such that
\be
\l\Vert \l[{\cal \hat M}_u^\a\r]^{-1}\r\Vert\le C^\a(u).\nn
\ee
\ele
\bpr
Let 
$
\l| (
\v,
\w
)
\r|_{H^{1,\a}(\R^2)\oplus L^{2,\a}(\R^2)} \le 1.
$
Due to \cref{co ds Rtwo alpha} there exists 
$
(
\t,
\p,
\lambda) \in \Big [ H^{2,\a}(\R^2)\oplus H^{1,\a}(\R^2) \Big ] \cap \ker {\cal L}_u^{\perp, L^{2,\a} \oplus L^{2,\a}} \oplus H^{2,\a}(\R)\,
$
such that
$
\bma
\v \\
\w 
\ema
= {\cal \hat M}_u^\a
(
\t,
\p ,
\lambda
)
$.
\cref{le Muinvestimate alpha} and the expression for
$(
\t,
\p,
\lambda
)$ from the proof of \cref{le invertibilityMx alpha} imply the claim.
\epr
\bde We define the following operators.
\begin{itemize}
\item [(a)] ${\cal K}_u^\a:  H^{2,\a}(\R^2)\oplus H^{1,\a}(\R^2)  \to H^{1,\a}(\R^2)\oplus L^{2,\a}(\R^2)$ given by \\
$$
\l({\cal K}_u^\a
\tp\r) (\xi,x) =
\bma
\u\dxi\t(\xi,x)-\p(\xi,x) \\ 
-\dx^2\t(\xi,x)+\cos(\t_K(\Z))\t(\xi,x)+\u\dxi\p(\xi,x) \\ 
\ema.
$$
\item [(b)] ${\cal N}_u^\a:  H^{2,\a}(\R^2)\oplus H^{1,\a}(\R^2) \oplus H^{2,\a}(\R)\to H^{1,\a}(\R^2)\oplus L^{2,\a}(\R^2) $ given by\\
$$
\l({\cal N}_u^\a
( \t ,
\p, 
\lambda )\r)(\xi,x)  =
 \l({\cal K}_{\u}^\a\tp\r) (\xi,x)
+ \lambda(\xi)
t_2(\xi,\u,x).$$
\item [(c)] 
${\cal \hat  N}_u^\a={\cal N}_u^\a\Bigg|_{\Big [ H^{2,\a}(\R^2)\oplus H^{1,\a}(\R^2) \Big ] \cap \l(\ker {\cal L}_u^\a\r)^{\perp, L^{2,\a} \oplus L^{2,\a}} \oplus H^{2,\a}(\R)} $.
\end{itemize}
\ede
\noindent Using von Neumann's theorem we are able to prove now that for small $u_*$ an extension of the operator ${\frak M}_n^\a:  Y_n^\a(u_*)  \to Z_n^\a(u_*)$ is invertible. Before proceeding to the proof we introduce the following definition in order to specify $u_*$.
\bde
Let $C^\a$ be a specific fixed function from \cref{le calMuinvestimate alpha}. Set
\be \la{ustarlow alpha}
{\ti u}^\a =\ti u\l(\l\Vert [\hat M^\a]^{-1}\r\Vert\r)=\sup\{u\in(-1,1)~|~\forall s,t\in \R: ~|s|,|t|\le |u| :|s|C^\a(t)<1 \}.
\ee
\ede
\bco\la{co invertibilityMxi alpha}
The linear operator
$
{\frak M}^\a:  Y^\a(u_*)  \to Z^\a(u_*),~
( \t,
\p,
\lambda) \mapsto
{\frak M}^\a
( \t,
\p,
\lambda),
$
given by
\be
{\frak M}^\a
( \t,
\p,
\lambda)(\u)=
 \hat {\cal  N}_u^\a ( \t(\u),
\t(\u),
\lambda (\u)),\nn
\ee
is invertible if $u_* < \ti  u^\a $.
\eco
\bpr
The operator $\ti{\frak M}^\a$
is invertible by \cref{le invertibilityMx alpha} 
and it holds for its operator norm that
$
\l\Vert \l[ \ti{\frak M}^\a\r]^{-1} \r\Vert \le \sup_{|\u|\le u_*}  C^\a(u)\,.
$
Let ${\frak P}^\a:  Y^\a(u_*)  \to Z^\a(u_*)$ be given by 
\be
{\frak P}^\a (\t,
\p,
\lambda) (u)=
\u\bma
\dxi\t(u)\\ 
\dxi\p(u) \\ 
\ema
+
\u\bma
\dx\t(u)\\ 
\dx\p(u) \\ 
\ema.\nn
\ee
It holds that 
$
\l\Vert {\frak P}^\a \r\Vert \le \sup_{|\u|\le u_*} |\u|\,
$
and thus
\be
\l\Vert {\frak P}^\a \r\Vert\l \Vert \l[\ti{\frak M}^\a\r]^{-1} \r\Vert 
\le  {}&  \sup_{|\u|\le u_*} |\u| \sup_{|\u|\le u_*}  C^\a(u)
<
1,\nn
\ee 
due to \re{ustarlow alpha}, since $u_* < \ti u^\a$.
Hence
${\frak P}^\a+\ti{\frak M}^\a={\frak M}^\a$ is invertible by von Neumann's theorem.
\epr
\subsection*{Completion of the Proof of \cref{le invertibilityMxiCtwo alpha}} 
Analogously to \cref{co invertibilityMxi alpha}
one shows first that the corresponding operator on spaces of higher regularity in $(\xi,x)$ is invertible. 
The claim  of \cref{le invertibilityMxiCtwo alpha} for the operator on spaces of higher regularity in $u$ and in $(\xi,x)$ follows by using difference quotients, orthogonal projection and the inverse mapping theorem. 
\qed
\paragraph{Acknowledgements}
My sincere gratitude goes to my PhD advisor  Markus Kunze for the continuous support, his patience and motivation.
I would like to express my deep appreciation to Justin Holmer for numerous helpful and fruitful discussions. 
%

\end{document}